\documentclass[fractalfract,article,accept,pdftex, moreauthors]{Definitions/mdpi} 
\firstpage{1} 
\pubvolume{1}
\issuenum{1}
\articlenumber{0}
\pubyear{2022}
\copyrightyear{2022}
\externaleditor{}
\datereceived{3 August 2022} 
\dateaccepted{13 September 2022} 
\datepublished{}

\hreflink{https://doi.org/} % If needed use \linebreak
\usepackage{float}
\usepackage[labelformat=simple]{subfig}

\usepackage{commath}
\usepackage{siunitx}
\usepackage{changes}
\newtheorem{thm}{Theorem}

\newtheorem{con}[thm]{Conjecture}
\newtheorem*{thm*}{Theorem}
\newtheorem*{lem*}{Lemma}
\newtheorem*{prop*}{Proposition}
\newtheorem*{con*}{Conjecture}

\theoremstyle{definition}

\newcommand*{\Com}{\mathbb{C}}
\newcommand*{\D}{\mathbb{D}}
\newcommand*{\M}{\mathcal{M}}

\newcommand*{\rs}{\hat{\mathbb{C}}}
\newcommand*{\A}{\mathcal{A}}

\newcommand*{\ra}{\rightarrow}
\newcommand*{\lra}{\longrightarrow}
\newcommand*{\sm}{\setminus}
\newcommand*{\ol}{\overline}

\Title{On Benford's Law and the Coefficients of the Riemann Mapping Function for the Exterior of the Mandelbrot Set}

\TitleCitation{On Benford's Law and the Coefficients of the Riemann Mapping Function for the Exterior of the Mandelbrot Set}

\Author{{Filippo Beretta} $^{1}$, Jesse Dimino $^{2,}$*, Weike Fang $^{3}$, Thomas C. Martinez $^{4}$, Steven J. Miller $^{5,}$* and Daniel Stoll $^{6}$} %MDPI: Please carefully check the accuracy of names and affiliations.

\AuthorNames{Filippo Beretta, Jesse Dimino, Weike Fang, Thomas C. Martinez, Steven J. Miller, and Daniel Stoll}

\AuthorCitation{Beretta, F.; Dimino, J.; Fang, W.; Martinez, T. C.; Miller, S. J.; Stoll, D.}
\address{%
$^{1}$ \quad Department of Mathematics, University of Milan, 20133 {Milan}, Italy; filo0820@gmail.com\\
$^{2}$ \quad Department of Mathematics, CUNY College of Staten Island, Staten Island, NY 10314, USA\\
$^{3}$ \quad Department of Mathematics, University of Notre Dame, {South Bend}, IN 46556, USA; wfang@nd.edu\\
$^{4}$ \quad Department of Mathematics, Harvey Mudd College, Claremont, CA 90701, USA; tmartinez@hmc.edu\\
$^{5}$ \quad Department of Mathematics and Statistics, Williams College, Williamstown, MA 01267, USA\\
$^{6}$ \quad Department of Mathematics, University of Michigan, Ann Arbor, MI 48109, USA; dastoll@umich.edu}

\corres{Correspondence: jesse.dimino@macaulay.cuny.edu (J.D.); sjm1@williams.edu (S.J.M.)}

\abstract{We investigate Benford's law in relation to fractal geometry. Basic fractals, such as the Cantor set and Sierpinski triangle are obtained as the limit of iterative sets, and the unique measures of their components follow a geometric distribution, which is Benford in most bases. Building on this intuition, we aim to study this distribution in more complicated fractals. We examine the Laurent {coefficients of a Riemann mapping} and {the} Taylor coefficients {of its} reciprocal function from the exterior of the Mandelbrot set to the complement of the unit disk. These coefficients are {2-adic} rational numbers, and through statistical testing, we demonstrate that the numerators and denominators are a good fit for Benford's law . We offer additional conjectures and observations about these coefficients. In particular, we highlight certain arithmetic subsequences related to the coefficients' denominators, provide an estimate for their slope, and describe efficient methods to compute them. }


\begin{document}

%%%%%%%%%%%%%%%%%%%%%%%%%%%%%%%%%%%%%%%%%%

% The order of the section titles is: Introduction, Materials and Methods, Results, Discussion, Conclusions for these journals: aerospace,algorithms,antibodies,antioxidants,atmosphere,axioms,biomedicines,carbon,crystals,designs,diagnostics,environments,fermentation,fluids,forests,fractalfract,informatics,information,inventions,jfmk,jrfm,lubricants,neonatalscreening,neuroglia,particles,pharmaceutics,polymers,processes,technologies,viruses,vision

\section{Introduction}
The Mandelbrot set $\M$ was first introduced and drawn by Brooks and Matelski. By~analyzing the family of functions $f_c(z) =  z^2 + c$, Douady and Hubbard began the formal mathematical study of the Mandelbrot set as the set of parameters $c$, for which the orbit of 0 under $f_c$ remains bounded. We study Benford's law in relation to the Mandelbrot set to both investigate the distribution's extension to fractal geometry and search for patterns in the Mandelbrot~set. 

In 1980, Douady and Hubbard were able to prove the connectedness of $\M$ by constructing a conformal isomorphism
\[\Phi: \Com\sm\M\lra\Com\sm\ol{\D}\]
between the complement of the Mandelbrot set and the complement of the{ closed unit disk}~\cite{DouH1}. Using the Douady-Hubbard map $\Phi$, we can define related conformal isomorphisms,
\begin{align*}
    &\Psi: {\Com}\sm\overline{\D} \lra {\Com}\sm\M,\\
    &\Theta: \D\lra\Com\sm\M^{-1},
\end{align*}
where $\M^{-1} =  \set{1/c: c\in\M}$, by~setting $\Psi =  \Phi^{-1}$ and $\Theta(c) =  1/\Psi(1/c)$.
One of the most heavily studied questions in complex dynamics is whether or not $\M$ is locally connected (MLC). By~a theorem of {Caratheodory}~\cite{LuoY1}, these two maps can be extended continuously to the unit circle if and only if the Mandelbrot set is locally connected. As~such, we focus on studying these maps and their respective Laurent and Taylor expansions, {as outlined in~\cite{Shi3} and~\cite{Shi2}, respectively}:
\begin{align}
   \Psi(z) &\ = \  z + \sum_{m=0}^\infty b_mz^{-m},\\
    \Theta(z) &\ = \  z + \sum_{m=2}^\infty a_mz^m.
\end{align}

In Section~\ref{alg_complx}, we outline the methods~\cite{Shi3,Shi2,EwiS2} we used to compute the $a_m$ and $b_m$ coefficients. The~computation time grows exponentially, so methods of improving computation are explored. Using recursion, we were able to compute the first 10240~coefficients.

In Section~\ref{benford_stats}, we discuss Benford's law along with statistical testing to determine whether the coefficients obey a Benford distribution. Given a base $b\ge2$, a~data set is Benford base $b$ if the probability of observing a given leading digit, $d$, is $\log_b\left((d+1)/d\right)$ (see~\cite{Benf1,Mil1}). {{We can write any positive $x$ as $S_b(x) b^{k_b(x)}$, where $S_b(x) \in [1, b)$ is the significand and $k_b(x)$ is an integer. If~the probability of observing a significand of $s \in [1, b)$ is $\log_b(s)$, we say the set is strongly Benford (or frequently just Benford).}} {In most cases, a~data set is demonstrated to be Benford through statistical testing. There are few straightforward proofs for Benfordness, all of which rely heavily on understanding the structure and properties of the data. Well understood sets such as geometric series and the Fibonacci numbers have explicit proofs for Benfordness, but~the structure and properties of coefficients we study are still the subject of active research in complex dynamics. Therefore, we rely on statistical testing for our results.} We consider the standard $\chi^2$ distribution and the sequence of the data's logarithms modulo 1 for our statistical testing, and~a standard goodness of fit test demonstrates that the numerators and denominators are a good fit for Benford's law, while the decimal representations are~not. 

\textls[-15]{Section~\ref{sec_coeffs} {deals with conjecture, observations, and~theorems related to the coefficients. Sections~\ref{subsec:amRemarks} and \ref{subsec:bmRemarks} are meant to tie together the most important of these that we have found from various authors for the $a_m$ and $b_m$ coefficients, respectively. In~Section~\ref{subsec:ourRemarks}}, we present new results and conjectures on the $a_m$ and $b_m$ coefficients {from our work}. Theorem~\ref{padic} \cite{Shi3,Shi4} {states} that they are 2-adic rational numbers; in other words, they are of the form $p/2^{-\nu}$, where $p$ is an odd integer. The~integer $\nu$ is, by~definition, the~$2$-adic valuation $\nu_2$ of $a_{m}$ or $b_{m}$.  Therefore, we focus on the denominator's exponents $-\nu(a_m)$, $-\nu(b_m)$. Setting $m = 2^nm_0$, with~$m_0$ odd and $n$ fixed; the subsequences $\{-\nu(a_m)\}$, $\{-\nu(b_m)\}$ appear to be near-arithmetic progressions. We present the results observed in the following~conjecture.}

\begin{con}\label{slope}
Let m be written as $2^nm_0$ as above, with~$n=\overline{n}$ fixed. Then, the sequence $\{-\nu(a_m)\}_{n=\overline{n}}$ is asymptotically linear, with~slope
\begin{align}
2/(2^{\overline{n}+1}-1).
\end{align}
\end{con}

We also present an efficient way to compute the denominator's exponents for the cases $n=0,1,2$. 

{Our work offers a new approach to some classical problems in complex dynamics, and~we do this through our statistical testing. These coefficients have been studied extensively, so we compile relevant observations from disparate authors to use as a basis for offering new conjectures and results. The~idea of studying Benford’s law in complex dynamics and fractal sets is also original; to our knowledge, no other authors have attempted to study Benford’s law in this setting. Our approach offers both a new and unique way to study important series in complex dynamics, and~it provides motivation for number theorists and statisticians to study Benford's law in the new field of data, namely fractal~sets.}

%%%%%%%%%%%%%%%%%%%%%%%%%%%%%%%%%%%%%%%%%%%%%%%%%%%%%%%%%%%%%%%%%%%%%%%%%%%%%%%%%%%%%%%%%%%%%%%%%%%%%%%%%%%%%%%%%%%%%%%%%%%%%%%%%%%%%%%%%%%%%%%%%%%%%%%%%%%%%%%%%%%%%%
%%%%%%%%%%%%%%%%%%%%%%%%%%%%%%%%%%%%%%%%%%%%%%%%%%%%%%%%%%%%%%%%%%%%%%%%%%%%%%%%%%%%%%%%%%%%%%%%%%%%%%%%%%%%%%%%%%%%%%%%%%%%%%%%%%%%%%%%%%%%%%%%%%%%%%%%%%%%%%%%%%%%%%
%%%%%%%%%%%%%%%%%%%%%%%%%%%%%%%%%%%%%%%%%%%%%%%%%%%%%%%%%%%%%%%%%%%%%%%%%%%%%%%%%%%%%%%%%%%%%%%%%%%%%%%%%%%%%%%%%%%%%%%%%%%%%%%%%%%%%%%%%%%%%%%%%%%%%%%%%%%%%%%%%%%%%%
\section{Preliminaries in Complex~Dynamics}

We give a brief introduction to complex dynamics. For~more detailed proofs see~\cite{Miln1} and~\cite{Sto1}. 

Let \(f:\rs\ra\rs\) be a rational map.
The \emph{Julia set} \(J_f\) associated to the map \(f\) may be defined as the closure of the set of repelling periodic points of \(f\).
For a rational map of degree 2 or higher, the~Julia set \(J_f\) is non-empty. We now restrict our attention to polynomial maps of degree \(d\ge 2\), which have a superattracting fixed point at infinity. We can thus define the \emph{filled Julia set} \(K_f\) as the complement of the basin of the attraction of infinity:
\[
K_f \ = \ \set{ z \in \Com: \set{f^n(z)}_{n\ge 0} \text{ is bounded}}.
\]

Making use of the above, it is possible to redefine \(J_f\) as the boundary of the filled Julia~set.

\begin{Lemma}
{Let f be a }polynomial of degree \(d \ge 2\). The~filled Julia set \(K_f \subset \Com\) is compact, with~boundary \(\partial K_f = J_f\) equal to the Julia set. The~complement \(\hat{\Com} \sm K_f\) is connected and equal to the basin of attraction \(A(\infty)\) of the point \(\infty\).
\end{Lemma}

It follows that the Julia set \(J_f\) of a polynomial \(f\) is precisely the boundary of the basin of attraction \(\A_f(\infty)\).

We may now characterize the connectedness of \(J_f\).
This is determined entirely by the activity of the critical points of \(f\).

\begin{Theorem}\label{Julia and Fil Julia}
 {The Julia set} \(J_f\) for a polynomial of degree \(d\ge2\) is connected if and only if the filled Julia set \(K_f\) contains every critical point of \(f\).
\end{Theorem}

%%%%%%%%%%%%%%%%%%%%%%%%%%%%%%%%%%%%%%%%%%%%%%%%%%%%%%%%%%%%%%%%%%%%%%%%%%%%%%%%%%%%%%%
%%%%%%%%%%%%%%%%%%%%%%%%%%%%%%%%%%%%%%%%%%%%%%%%%%%%%%%%%%%%%%%%%%%%%%%%%%%%%%%%%%%%%%%
%%%%%%%%%%%%%%%%%%%%%%%%%%%%%%%%%%%%%%%%%%%%%%%%%%%%%%%%%%%%%%%%%%%%%%%%%%%%%%%%%%%%%%%
%%%%%%%%%%%%%%%%%%%%%%%%%%%%%%%%%%%%%%%%%%%%%%%%%%%%%%%%%%%%%%%%%%%%%%%%%%%%%%%%%%%%%%%
\subsection*{The Mandelbrot~Set}\text{}

We now focus primarily on the family of quadratic functions of the form \(\{f_c(z) = z^2 + c\}_{c \in \mathbb{C}}\). 
Since \(f_c\) has a single critical point 0, it follows from {Theorem~\ref{Julia and Fil Julia}} that \(J_{f_c}\) is connected if and only if the orbit \(\set{f_c^n(0) \mid n\in \mathbb{N}}\) is bounded. This motivates the following definition of the Mandelbrot set, \(\M\).

\begin{Definition}
\(\M\subset\Com\){ is the set of all the parameters} \(c\in\Com\) such that the Julia set \(J_{f_c}\) is connected. 
Equivalently, \(\M\) is the set of all the \(c\) such that the orbit of 0 under \(f_c\) remains bounded:
\[
    \M \ := \ \set{c\in\Com:\,\text{there exists } R>0\text{ such that for all }n,|f^n_c(0)|<R}.
\]
\end{Definition}

\begin{Remark}
{It is possible to generalize this} definition and most of the following results to the family of unicritical degree $d$ polynomials \(f_{c,d}(z) = z^d + c\), where \(d\) is an integer \(d \ge 2\).
In this case, \(\M_d\) is called the \emph{multibrot set} of degree \(d\).
For simplicity, we focus only on \(\M = \M_2\), which has historically been the object of greatest interest.
For more information on \(\M_d\), see~\cite{Shi3,Shi4}.
\end{Remark}

\begin{figure}[H]
    
    \includegraphics[height = 8cm, trim={0 5cm 0 5cm}, clip] {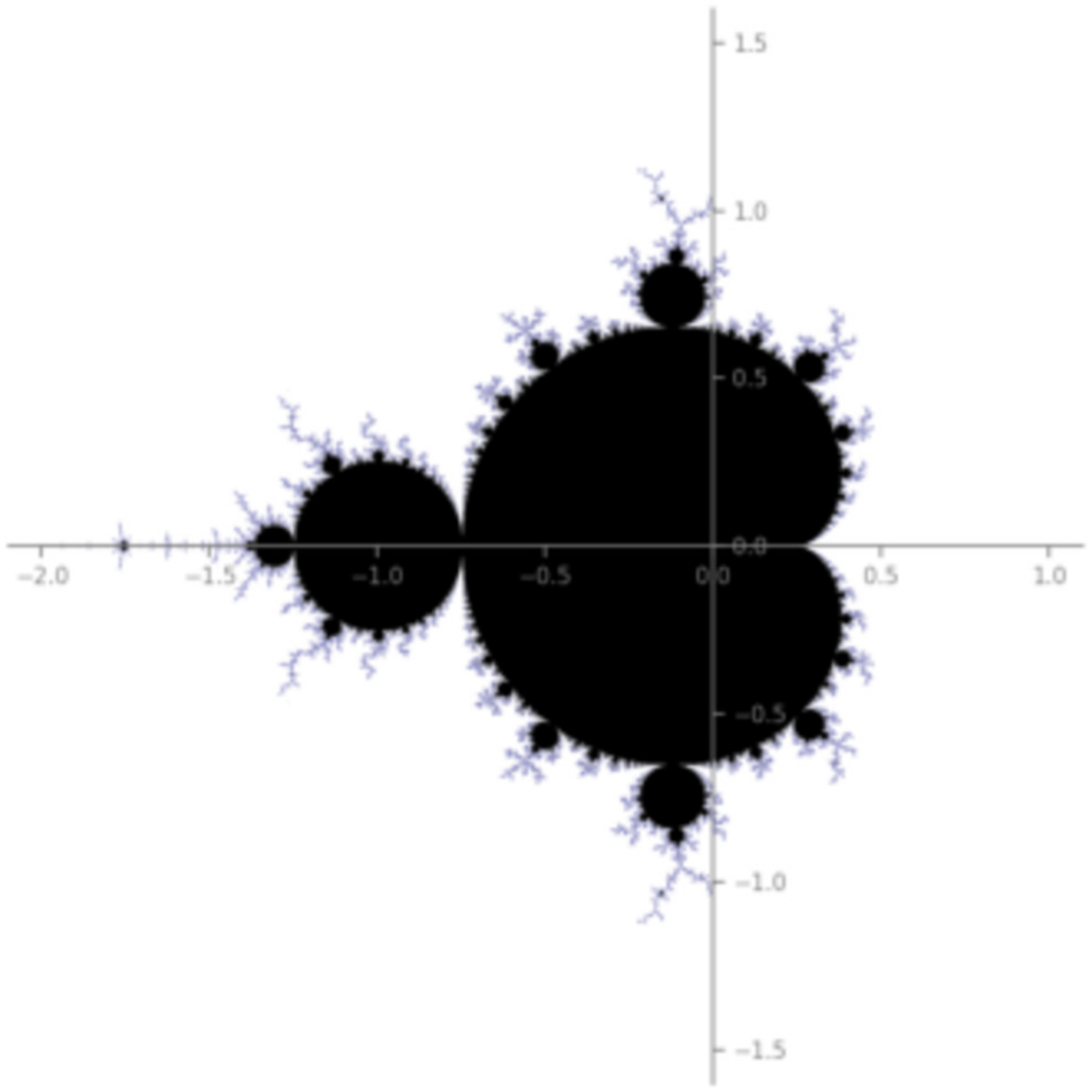}
    \caption{The Mandelbrot set \(\M\) in the complex~plane.}
    \label{fig:MandelbrotSet}

\end{figure}

It is possible to demonstrate that the interior of $\mathcal{M}$ is nonempty. We utilized an escape time algorithm and computer graphics to obtain the visualization of \(\M\) presented in Figure \ref{fig:MandelbrotSet}. 

When the first computer images of \(\M\) were generated, Benoit Mandelbrot observed small regions that appeared to be separate from the main cardioid, and~conjectured that \(\M\) was disconnected, which was later~disproved.

\begin{Theorem}[{Douady, Hubbard}]
 The Mandelbrot set is connected.
\end{Theorem}

This result was first proved by Douady and Hubbard~\cite{DouH1} by explicitly constructing a conformal isomorphism \(\Phi: \rs\sm\M \to \rs\sm\ol{\D}\).
Douady and Hubbard's proof is significant not only for the result, but~also since it provides an explicit formula for the uniformization of the complement of the Mandelbrot~set.

A large amount of research has been devoted to the local connectivity of the Mandelbrot set, which is generally regarded as one of the most important open problems in complex dynamics.
We recall that a set \(A\) in a topological space \(X\) is locally connected at \(p \in A\) if for every open set \(V\subset X\) containing \(p\), there is an open subset \(U\) with \(p \in U \subset V\) such that \(U \cap X\) is connected.
The set $A$ is said to be \emph{locally connected} if it is locally connected at $p$ for all \(p \in A\).

As above, let \(\Phi: \hat{\Com} \sm \M \to \hat{\Com} \sm \ol{\D}\) be the conformal isomorphism constructed by Douady and Hubbard; notice that the map \(\Psi: \hat{\Com}\sm\overline{\D} \to \hat{\Com}\sm\M\) is the Riemann mapping function of \(\hat{\Com}\sm\M\).
%, as defined in \emph{Theorem \ref{Rie Map}}. 
We consider it the Laurent expansion at \(\infty\) 
\[\Psi(z) \ = \ z +\sum_{m=0}^\infty b_mz^{-m}.\]

Another possibility is to consider \(\Theta(z) :=1/\Psi(1/z)\), which is the Riemann mapping of the bounded domain \(\Com\sm\{1/z : z \in \M\}\). We have the corresponding Taylor expansion for $\Theta$ at the origin:
\[\Theta(z) \ = \ z+\sum_{m=2}^\infty a_mz^m.\]

{For the general Multibrot set $\mathcal{M}_d$, we refer to the coefficients with the notation $b_{d,m}$ and $a_{d,m}$.} To underline the importance of these maps, we reference a lemma from Caratheodory~\cite{LuoY1}.

\begin{Theorem} 
\label{car}
 \textbf{(Caratheodory's continuity lemma)} Let $G\subset\hat{\Com}$ be a simply connected domain and a function $f$ maps $\D$ onto $G$ via a conformal isomorphism. Then $f$ has a continuous extension to $\overline{\D}$ if and only if the boundary of $G$ is locally connected.
\end{Theorem}

Therefore, if~the map $\Psi$ or the map $f$ can be extended continuously to the unit circle $\partial \D$, then $\M$ is locally connected. To~demonstrate this extension, it would be sufficient to prove that one of the two series converges uniformly on $\D$.

There have been numerous attempts to prove this result. For~example, Ewing and Schober demonstrated that the inequality $0 < |b_m|$ < 1/m holds for every m < $240,000$. A~bound of the type $|b_m| < K/m^{(1+\epsilon)}$ would lead to the desired result; however, this would imply that the extension of $\Psi$ is Hölder continuous, which it is not, as expressed in~\cite{BieFH1}. Proving that $|b_m| <$ K/(m$\;\log^2$(m))  would prove that the series converges absolutely; however, modern computations suggest that such a bound does not exist~\cite{BieFH1}. Therefore, the~MLC conjecture and its consequences remain an object of active study.

\begin{Remark}[Consequences of local connectedness]
This paper tackles problems related {to the} behaviour of the $\{a_m\}$ and $\{b_m\}$ coefficients and the local connectedness of $\M$. In~this section, we highlight the importance of this property of $\M$. In~particular, it implies a conjecture known as the Density of Hyperbolicity. A~\emph{hyperbolic component} is an interior connected component of $\M$, in which the sequences $\{f_c^n(0)\}$ have an attracting periodic cycle of period $p$. The~\emph{Density of Hyperbolicity} conjecture states that these are the only interior regions of $\M$. Another consequence of MLC is related to a topologically equivalent description of $\partial\M$. In~particular, the~boundary of the Mandelbrot set can be identified with the unit circle $S_1$ under a specific relation $\thicksim$, known as the \emph{abstract Mandelbrot set}~\cite{BanK1}. More information and other implications of MLC and the Density of Hyperbolicity may be found in~\cite{Beni1} and~\cite{DouH1}.
\end{Remark}

%%%%%%%%%%%%%%%%%%%%%%%%%%%%%%%%%%%%%%%%%%%%%%%%%%%%%%%%%%%%%%%%%%%%%%%%%%%%%%%%%%%%%%%%%%%%%%%%%%%%%%%%%%%%%%%%%%%%%%%%%%%%%%%%%%%%%%%%%%%%%%%%%%%%%%%%%%%%%%%%%%%%%%
%%%%%%%%%%%%%%%%%%%%%%%%%%%%%%%%%%%%%%%%%%%%%%%%%%%%%%%%%%%%%%%%%%%%%%%%%%%%%%%%%%%%%%%%%%%%%%%%%%%%%%%%%%%%%%%%%%%%%%%%%%%%%%%%%%%%%%%%%%%%%%%%%%%%%%%%%%%%%%%%%%%%%%
%%%%%%%%%%%%%%%%%%%%%%%%%%%%%%%%%%%%%%%%%%%%%%%%%%%%%%%%%%%%%%%%%%%%%%%%%%%%%%%%%%%%%%%%%%%%%%%%%%%%%%%%%%%%%%%%%%%%%%%%%%%%%%%%%%%%%%%%%%%%%%%%%%%%%%%%%%%%%%%%%%%%%%
\section{Algorithms and~Complexity}\label{alg_complx}

There are algorithms to compute both the $a_m$ and $b_m$ coeffecients. While these algorithms work for a generic degree $d\ge2$, we focus on $d=2$, which is {historically} the most interesting case, since it is the one associated with $\mathcal{M}$. For~simplicity, we denote $b_m = b_{2,m}$ and $a_m = a_{2,m}$. The~behavior for other values is similar. Derivations for the explicit form of the $b_m$  and $a_m $ coefficients may be found in~\cite{Shi3,Shi2} respectively. 
There is also a formula to switch between the coefficients, outlined in~\cite{Shi1}. 

\begin{Theorem} \label{alg1}
    Let $n\in\mathbb{N}$ and $1\le m \le d^{n+1}-3$. Then 
    \begin{align*} 
    \centering
    b_{d,m} \ = \  - \frac{1}{m} \sum &C_{j_1}\left(\frac{m}{d^n}\right) C_{j_2}\left(\frac{m}{d^{n-1}} - dj_1\right)C_{j_3}\left(\frac{m}{d^{n-2}} - d^2j_1 - dj_2\right)\cdots\\
    &C_{j_n}\left(\frac{m}{d} - d^{n-1}j_1-\cdots - dj_{n-1}\right),
    \end{align*} %\fix{FIXED - MILLER: DO WE NEED A MINUS SIGN AFTER $j_1$ AND BEFORE THE $\cdots$?}
    where the sum is over all non-negative indices $j_1,\dots, j_n$ such that
\begin{equation}
        \left(d^n - 1\right)j_1 +
    \left(d^{n-1}-1\right)j_2 +\left(d^{n-2} -1\right)j_3 +\cdots+\left(d-1\right)j_n \ = \ m+1
    \label{dio1}
    \end{equation}
    and $C_j\left(a\right)$ is the generalized binomial coefficient
    \begin{align*}
        C_j\left(a\right) \ = \  \frac{a\left(a-1\right)\cdots\left(a-j+1\right)}{j\left(j-1\right)\cdots 2 \cdot 1}.
    \end{align*}
\end{Theorem} 

The $a_m$ coefficients can be obtained from the $b_m$ using the formula
\begin{equation}
\label{abrelation}
    a_{d,m}\ = \ -b_{d,m-2}- \sum_{j=2}^{m-1}a_{d,j}b_{d,m-1-j},
\end{equation}
or they can be directly calculated as in the following~theorem.

\begin{Theorem} \label{alg2}
    Let $n\in\mathbb{N}$ and $2\le m \le d^{n+1}-1$. Then 
    \begin{align*}
    \centering
    a_{d,m} \ = \  \frac{1}{m} \sum &C_{j_1}\left(\frac{m}{d^n}\right) C_{j_2}\left(\frac{m}{d^{n-1}} - dj_1\right)C_{j_3}\left(\frac{m}{d^{n-2}} - d^2j_1 - dj_2\right)\cdots\\
    &C_{j_n}\left(\frac{m}{d} - d^{n-1}j_1 -  \cdots - dj_{n-1}\right),
    \end{align*}
    where the sum is over all non-negative indices $j_1,\dots, j_n$ such that
\begin{equation}
    \label{dio2}
        (d^n - 1)j_1 +
    (d^{n-1}-1)j_2 +(d^{n-2} -1)j_3 +\cdots+(d-1)j_n \ = \ m - 1
    \end{equation}
    and $C_j(a)$ is the generalized binomial coefficient
    \begin{align*}
        C_j(a) \ = \  \frac{a(a-1)\cdots(a-j+1)}{j(j-1)\cdots2 \cdot 1}.
    \end{align*}
\end{Theorem}

While the above theorems give the explicit forms of the coefficients, the~following theorem provides a recursive method to find $b_m$, which is more suitable for computers. Once we find $b_m$, we can apply the relationship between $a_m$ and $b_m$ outlined in theorem~\ref{alg1} to find $a_m$. More details can be found in~\cite{EwiS2}.

\begin{Theorem} \label{alg3}
    Let $n\in\mathbb{N}$ and $d=2$.
    Then
    \begin{align*}
    \centering
    b_{m} \ = \ \beta_{0,m+1}
    \end{align*}
    where the following holds true.
\[\beta_{j,1} = 0, j \geq 1; \beta_{n,0} = 1, n \in \mathbb{N}\]
\[\beta_{n,m} = 0, n \geq 1, 1\leq m \leq 2^{n+1} - 2\]
    \begin{align*}
    \centering
    \beta_{n-1, m} \ = \ \frac{1}{2} \left[ \beta_{n,m} - \sum_{k = 2^n-1}^{m-2^n+1} \beta_{n-1,k} \beta_{n-1,m-k} - \beta_{0, m-2^n+1}  \right]
    \end{align*}

\end{Theorem}

\begin{Example}
    For example, to~calculate the first several $b_m$ coefficients, we can use Theorem~\ref{alg3} to obtain:
\[
    b_{0} \ = \ \beta_{0,1} \ = \ \frac{1}{2} \left[ \beta_{1,1} - \beta_{0, 0} \right] \ = \ \frac{1}{2} \left[ 0 - 1 \right] \ = \ -\frac{1}{2}
\]
\[
    b_{1} \ = \ \beta_{0,2} \ = \ \frac{1}{2} \left[ \beta_{1,2} - \beta_{0,1}^2 - \beta_{0, 1} \right] \ = \ \frac{1}{2} \left[ 0 - 1/4 + 1/2 \right] \ = \ \frac{1}{8}
\]
\[
    b_{2} \ = \ \beta_{0,3} \ = \ \frac{1}{2} \left[ \beta_{1,3} - 2(\beta_{0,1} \beta_{0,2}) - \beta_{0, 2} \right] \ = \ \frac{1}{2} \left[ \beta_{1,3} + 1/8 - 1/8 \right] = - \frac{1}{4}
\]
    since         
    \begin{align*}
    \centering
    \beta_{1,3} \ = \ \frac{1}{2} \left[ \beta_{2,3} - \beta_{0, 0} \right] \ = \ \frac{1}{2} \left[ 0 - 1 \right] \ = \ -\frac{1}{2}
    \end{align*}
\end{Example}

%%%%%%%%%%%%%%%%%%%%%%%%%%%%%%%%%%%%%%%%%%%%%%%%%%%%%%%%%%%%%%%%%%%%%%%%%%%%%%%%%%%%%%%
%%%%%%%%%%%%%%%%%%%%%%%%%%%%%%%%%%%%%%%%%%%%%%%%%%%%%%%%%%%%%%%%%%%%%%%%%%%%%%%%%%%%%%%
%%%%%%%%%%%%%%%%%%%%%%%%%%%%%%%%%%%%%%%%%%%%%%%%%%%%%%%%%%%%%%%%%%%%%%%%%%%%%%%%%%%%%%%
%%%%%%%%%%%%%%%%%%%%%%%%%%%%%%%%%%%%%%%%%%%%%%%%%%%%%%%%%%%%%%%%%%%%%%%%%%%%%%%%%%%%%%%
\subsection{Direct~Computation}

{Our initial algorithm to generate these coefficients was to directly compute them, and~we include our original methodology for reference as we cross checked our results. Most of the data was generated through the recursive algorithm, and~details are provided in Section~\ref{sec:Rec}.}

We wrote a program in Python to compute $b_{2,m}$ and $a_{2,m}$ based on the formulas given in Theorems~\ref{alg1} and ~\ref{alg2}, and~we obtained the first 1024 exact values of both coefficients'~sequences. 

Our methodology for computing the $m \textsuperscript{th}$ coefficient was to first generate the solutions $j_1, \dots, j_n$ to the Diophantine Equations \eqref{dio1} and \eqref{dio2}. Then, we plug them into the exact formula of $a_{m} = a_{2,m}$ and $b_{m} = b_{2,m}$ to find the $m \textsuperscript{th}$ coefficient.

We improve{d} the method to solve the Diophantine equations by first setting an upper limit on the degree for which to obtain coefficients, generate{d} the solutions for the highest order coefficient, and~create{d} data structures for dynamic storage. We stored each individual solution as a tuple of length $n$, where the $k^{th}$ entry denoted the value of the coefficient $j_k$. Every solution for the upper bound was then given a reference in a linked list, which we can use to find the highest order coefficient. The~solution stored for the upper bound can then be modified through decrementing the value for $j_n$ for each tuple and then deleting the reference to the tuple and deallocating the memory in the linked list when the value for $j_n$ reaches~zero. 

To deal with the time sink in generating the binomial coefficients, we utilized multi-core parallel computing. Each coefficient can be computed independently after we obtain the solutions to the Diophantine equations. We structure our code for concurrent computation and use generator expressions so that we can use multiple cores where our code is executing simultaneously. In~a multi-core setting, each core deals with one coefficient at a time. When one coefficient calculation is conducted, the~core takes the next awaiting task that is not being taken by other cores. We chose a high-performance server machine and ran our code in a parallel environment with $72$ valid cores. We obtained the first 1024 coefficients with a CPU time of $166$ hours and a total run time of $7$ days.

\subsection{Recursive~Computation}\label{sec:Rec}

The direct computation runs in exponential time, and~it is {generally impractical for generating large degree coefficients.} Therefore, we switched to a recursive method to generate these~coefficients. 

The method is described in~\cite{Shi1, BieFH1, EwiS2}, and~we outline the formula for the computation in Theorem~\ref{alg3}. {This method is efficient because it is able to reuse} information {from the previous coefficients} to compute the next {one}.

We wrote a Python program to implement the recursion to find $b_m$ and then use Equation \eqref{abrelation} to find the corresponding $a_m$. We were able to obtain {10240} coefficients for both series within {82} hours with a single core. We are also able to cross-check our computation results with the direct computation method before starting the statistical analysis. {All codes and results can be found at \url{https://github.com/DannyStoll1/polymath-fractal-geometry}  (accessed on 19 September 2022}). Detailed instructions can be found in the README file. 

%%%%%%%%%%%%%%%%%%%%%%%%%%%%%%%%%%%%%%%%%%%%%%%%%%%%%%%%%%%%%%%%%%%%%%%%%%%%%%%%%%%%%%%%%%%%%%%%%%%%%%%%%%%%%%%%%%%%%%%%%%%%%%%%%%%%%%%%%%%%%%%%%%%%%%%%%%%%%%%%%%%%%%
%%%%%%%%%%%%%%%%%%%%%%%%%%%%%%%%%%%%%%%%%%%%%%%%%%%%%%%%%%%%%%%%%%%%%%%%%%%%%%%%%%%%%%%%%%%%%%%%%%%%%%%%%%%%%%%%%%%%%%%%%%%%%%%%%%%%%%%%%%%%%%%%%%%%%%%%%%%%%%%%%%%%%%
%%%%%%%%%%%%%%%%%%%%%%%%%%%%%%%%%%%%%%%%%%%%%%%%%%%%%%%%%%%%%%%%%%%%%%%%%%%%%%%%%%%%%%%%%%%%%%%%%%%%%%%%%%%%%%%%%%%%%%%%%%%%%%%%%%%%%%%%%%%%%%%%%%%%%%%%%%%%%%%%%%%%%%
\section{Benford's~Law}\label{benford_stats}

\begin{figure}[H]
    
    \includegraphics[height = 6cm, width = 10cm]{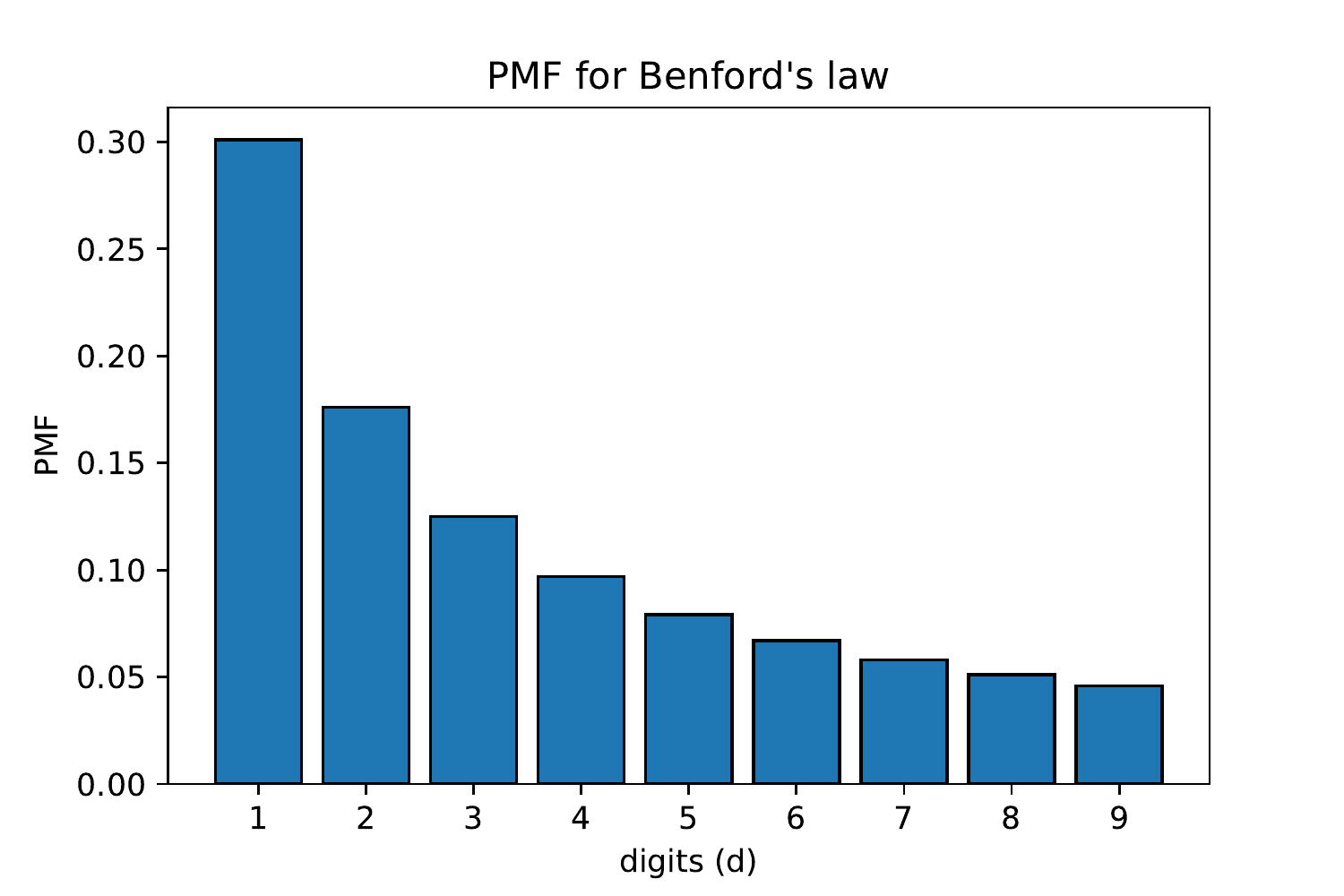}
    \caption{Probability Measure for Benford's Law in Base~10.}
    \label{fig:BenfordPM}
\end{figure}

Frank Benford's 1938 paper, \textit{The Law of Anomalous Numbers}~\cite{Benf1}, illustrated a profound result, in which the~first digits of numbers in a given data set are not uniformly distributed in general. Benford applied statistical analysis to a variety of well-behaved but uncorrelated data sets, such as the areas of rivers, financial returns, and~lists of physical constants; in an overwhelming amount of the data, 1 appeared as the leading digit around 30\% of the time, and~each higher digit was successively less likely~\cite{Benf1, Mil1}. He then outlined the derivation of a statistical distribution which maintained that the probability of observing a leading digit, $d$, for~a given base, $b$, is $\log_b\left((d+1)/d\right)$ for such data sets~\cite{Benf1, Mil1}. This logarithmic relation is referred to as Benford's law, and~its resultant probability measure for base 10 is outlined in Figure~\ref{fig:BenfordPM}. Benford's law has been the subject of intensive research over the past several decades, arising in numerous fields; see~\cite{Mil1} for an introduction to the general theory and numerous~applications. 

%\textcolor{red}{*Note for Miller to add more/better references}

Benford's law appears throughout purely mathematical constructions such as geometric series, recurrence relations, and~geometric Brownian motion. Its ubiquity makes it one of the most interesting objects in modern mathematics, as~it arises in many disciplines. Therefore, it {is} worthwhile to consider non-traditional data, such as fractals, where {the distribution} may appear. Basic fractals such as the Cantor set and Sierpinski triangle are obtained as the limits of iterations on sets, and their component measures (the lengths in the Cantor set and the areas in the Sierpinski triangle) follow a geometric distribution, which {is} Benford in most bases. Building on these results, it is plausible that more complicated fractals obey this distribution as well. We studied the Riemann mapping of the {exterior} of the Mandelbrot set to the {complement} of the unit disk, along with its reciprocal function to determine their fit to Benford's law. These mappings are given by a Taylor and Laurent series, respectively. The coefficients are of interest as their asymptotic convergence is intimately related to the conjectured local connectivity of the Mandelbrot set, which is an important open problem in complex~dynamics.

%%%%%%%%%%%%%%%%%%%%%%%%%%%%%%%%%%%%%%%%%%%%%%%%%%%%%%%%%%%%%%%%%%%%%%%%%%%%%%%%%%%%%%%
%%%%%%%%%%%%%%%%%%%%%%%%%%%%%%%%%%%%%%%%%%%%%%%%%%%%%%%%%%%%%%%%%%%%%%%%%%%%%%%%%%%%%%%
%%%%%%%%%%%%%%%%%%%%%%%%%%%%%%%%%%%%%%%%%%%%%%%%%%%%%%%%%%%%%%%%%%%%%%%%%%%%%%%%%%%%%%%
%%%%%%%%%%%%%%%%%%%%%%%%%%%%%%%%%%%%%%%%%%%%%%%%%%%%%%%%%%%%%%%%%%%%%%%%%%%%%%%%%%%%%%%
\subsection{Statistical Testing for Benford's~Law}

A common practice for evaluating whether a data set is distributed according to Benford's law is to utilize the standard $\chi^2$ goodness of fit test. As~we are investigating Benford's law in base 10, we utilize 8 degrees of freedom for our $\chi^2$ testing. ({{There are nine possible first digits, but~once we know the proportion that {is} digits 1 through 8 the percentage that starts with a 9 is forced, and~thus we lose one degree of freedom}}). If there are $N$ observations, letting $p_d = \log_{10} \left((d+1)/d\right)$ be the Benford probability of having a first digit of $d$, we expect $p_d N$ values to have a first digit of $d$. If~we let $\mathcal{O}_d$ be the observed number with a first digit of $d$, the~$\chi^2$ value is $$\chi^2 \ := \ \sum_{d=1}^9 \frac{(\mathcal{O}_d - p_d N)^2}{p_d N}. $${If the data are Benford, w}ith 8 degrees of freedom,  then 95\% of the time the $\chi^2$ test will produce a value of at most 15.5073; {this corresponds to a significance level of $\alpha = 0.05$. We perform multiple testing by creating a distribution of $\chi^2$ values as a function of sample size up to the m\textsuperscript{th} coefficient. This is standard practice for studying Benford sequences, and this is done to account for periodicity in the $\chi^2$ values, which is typical for certain Benford data sets, such as the integer powers of 2. To~account for this multiplicity, we also incorporate the standard Bonferroni correction. The overall testing is conducted at the level of significance of $\alpha = 0.05$, while giving equal weight in terms of significance to each individual test by conducting them at a significance level of $\alpha/m$. The~rationale is to keep the significance of the overall test constant with respect to the number of tests performed. As~we increase the total number of tests performed, we wish to increase our correction accordingly.} In total, we perform {10,045} tests for the $a_m$ and {10,046} tests for the $b_m$ coefficients as we compute the $\chi^2$ statistic each time we add a new non-zero coefficient to our data set. This brings our corrected threshold value to {38.9706 and 38.9708}, respectively. This corresponds to $\alpha$ = 0.05/10,045 = 4.978 $\times$ 10$^{-6}$ for each individual hypothesis in the $a_m$ dataset and $\alpha$ = 0.05/10,046 = 4.977$\times$10$^{-6}$ for each individual hypothesis in the $b_m$ dataset.

{Each data point is not independent, as the $\chi^2$ values are computed by using a running total of the data, and~as such, each point is built on the previous one. This results in a high correlation between the data, and~ the Bonferroni correction likely overcompensates for the increase in type I error. Still, it is one of the most plausible methods of dealing with the increase in multiplicity, since it is one of the simplest and most conservative estimates, and~a value above the Bonferroni correction provides strong evidence that the data are not Benford. Using independent increments to compute each $\chi^2$ statistic for Benford's law would fix the issue of independence, but~is not recommended since periodic behavior can be missed if the increments are chosen poorly.}

We considered the distribution of the $\chi^2$ values to account for random fluctuations and periodic {effects}. {In addition, we provide the provide the \emph{p}-value of our computed $\chi^2$ statistic to provide the type I error rate for our conclusions, and~we compute the powers of the $\chi^2$ tests relative to our null hypothesis that the data are Benford by using the noncentral chi-squared distribution~\cite{Guenther}.} {We also conducted simulations to estimate the sampling error relative to our null hypothesis. We wish to see how likely it is that a random sample falls in the rejection region for our testing; based on our significance level of $\alpha =0.05$, we expect the sampling error to be roughly 5\% if the data are Benford. We randomly sample 1000 coefficients from our data sets with replacement and calculate the $\chi^2$ value for this sample data. We repeat this simulation 1000 times, and~take the ratio the values in the rejection region to the total number of sample statistics calculated, to~estimate the sampling error.}

An equivalent test is to consider the distribution of the base 10 logarithm{s} {for} the {absolute value of the }data set modulo 1; a necessary and sufficient condition for a Benford distribution is that this sequence converges to a uniform distribution~\cite{MilTB1}. {To quantify the uniformity of this distribution, we again consider the standard $\chi^2$ test. We only perform a single test for the total data set, so we do not need to account for multiplicity.  Specifically, we split the interval $[0,1]$ into 10 equal bins. If~the data are uniform, we expect that each bin obtains $1/10$ of the total data. Therefore, for~each bin in the $a_m$ data we expect a value of $10,045/10 = 1004.5$ and for each bin in the $b_m$ data we expect a value of $10,046/10 = 1004.6$. Since there are 10 possible observations, we have 9 degrees of freedom for the data. If~the percentage of the data in the first nine bins is known, then the percentage of data in the last bin is forced, and~we lose one degree of freedom. We again take $\alpha=0.05$ and for nine degrees of freedom, this corresponds to a $\chi^2$ value of 16.919. These $\chi^2$ results are generated by cells 16, 17, and~18 by the Jupyter notebooks amLogData.ipynb and bmLogData.ipynb, respectively, which may be found under the Data Analysis folder. We also provide the associated \emph{p}-values and powers of the $\chi^2$ statistic relative to the null hypothesis that the data are uniform.}

The coefficients we studied are {2-adic }rational numbers, so we considered the distributions of the numerators, denominators, and~decimal expansions separately. We considered {only }the non-zero coefficients, since zero is not defined for our probability measure, and~certain theorems and conjectures outlined by Shiamuchi in~\cite{Shi2} already describe the distribution of the zeroes in the coefficients. {Our goal is to identify which components of these coefficients are a good fit for Benford's law through statistical testing.}

{Table \ref{Tab_coeff} provides} example{s} of the coefficients computed. When the coefficient is 0, numerators are set to 0 and denominators to - {for readability}. We then use them to compute the exact values in decimal expansion for $a_m$ and $b_m$.

\begin{table}[H]
\caption{{The table of firs}t {11} coefficients computed.{(codes/coeff\_compute.py)}.}
\label{Tab_coeff}
\newcolumntype{C}{>{\centering\arraybackslash}X}
\begin{tabularx}{\textwidth}{CCCCC}
\toprule
 \boldmath{$m$} & \boldmath{$a_m$} \textbf{Num} & \boldmath{$a_m$} \textbf{Denom} & \boldmath{$b_m$} \textbf{Num} & \boldmath{$b_m$} \textbf{Denom}\\ [0.5ex] 
 \midrule
0 & 0 & - & $-$1 & 2\\
 \midrule
 1 & 0 & -  & 1 & 8 \\ 
 \midrule
 2 & 1 & 2 & $-$1 & 4\\
 \midrule
 3 & 1 & 8 & 15 & 128\\
 \midrule
 4 & 1 & 4 & 0 & - \\
 \midrule
 5 & 15 & 128 & $-$47 & 1024\\
 \midrule
 6 & 0 & -  & $-$1 & 16\\ 
 \midrule
 7 & 81 & 1024 & 987 & 32,768 \\
 \midrule
 8 & 1 & 8 & 0 & -\\
 \midrule
 9 & 1499 & 32,768 & $-$3673 & 262,144\\
 \midrule
 10 & 1 & 32 & 1 & 32\\
 \bottomrule
\end{tabularx}
%\end{center}
%\vspace{15pt}
%
%\vspace{-15pt}
\end{table}

%%%%%%%%%%%%%%%%%%%%%%%%%%%%%%%%%%%%%%%%%%%%%%%%%%%%%%%%%%%%%%%%%%%%%%%%%%%%%%%%%%%%%%%
%%%%%%%%%%%%%%%%%%%%%%%%%%%%%%%%%%%%%%%%%%%%%%%%%%%%%%%%%%%%%%%%%%%%%%%%%%%%%%%%%%%%%%%
%%%%%%%%%%%%%%%%%%%%%%%%%%%%%%%%%%%%%%%%%%%%%%%%%%%%%%%%%%%%%%%%%%%%%%%%%%%%%%%%%%%%%%%
%%%%%%%%%%%%%%%%%%%%%%%%%%%%%%%%%%%%%%%%%%%%%%%%%%%%%%%%%%%%%%%%%%%%%%%%%%%%%%%%%%%%%%%
\subsection{Benfordness of the Taylor and Laurent~Coefficients}

\begin{figure}[H]
    \subfloat[\centering $a_m$]{{\includegraphics[height = 6cm, width = 6cm]{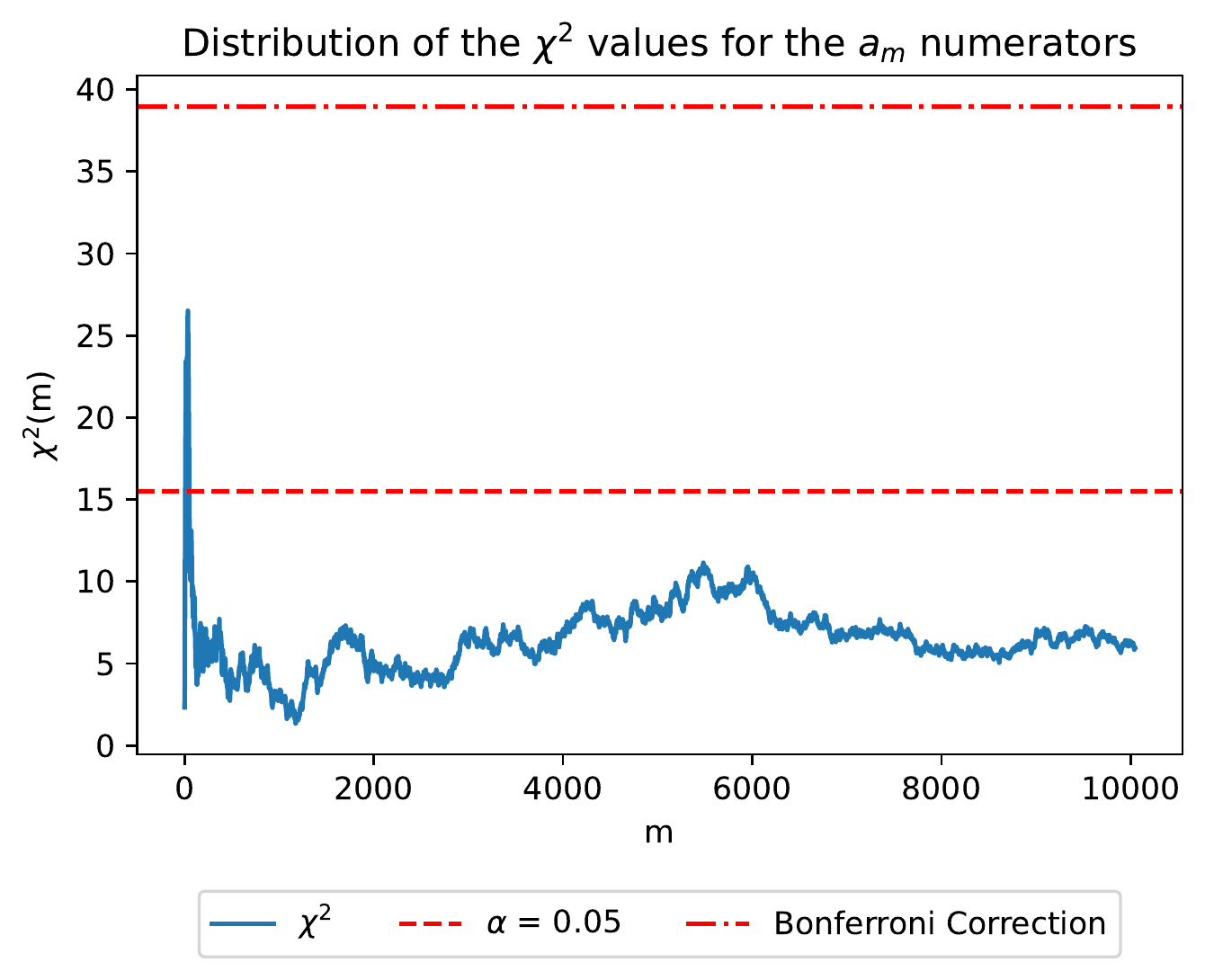} }}%
    \qquad
    \subfloat[\centering $b_m$]{{\includegraphics[height = 6cm, width = 6cm]{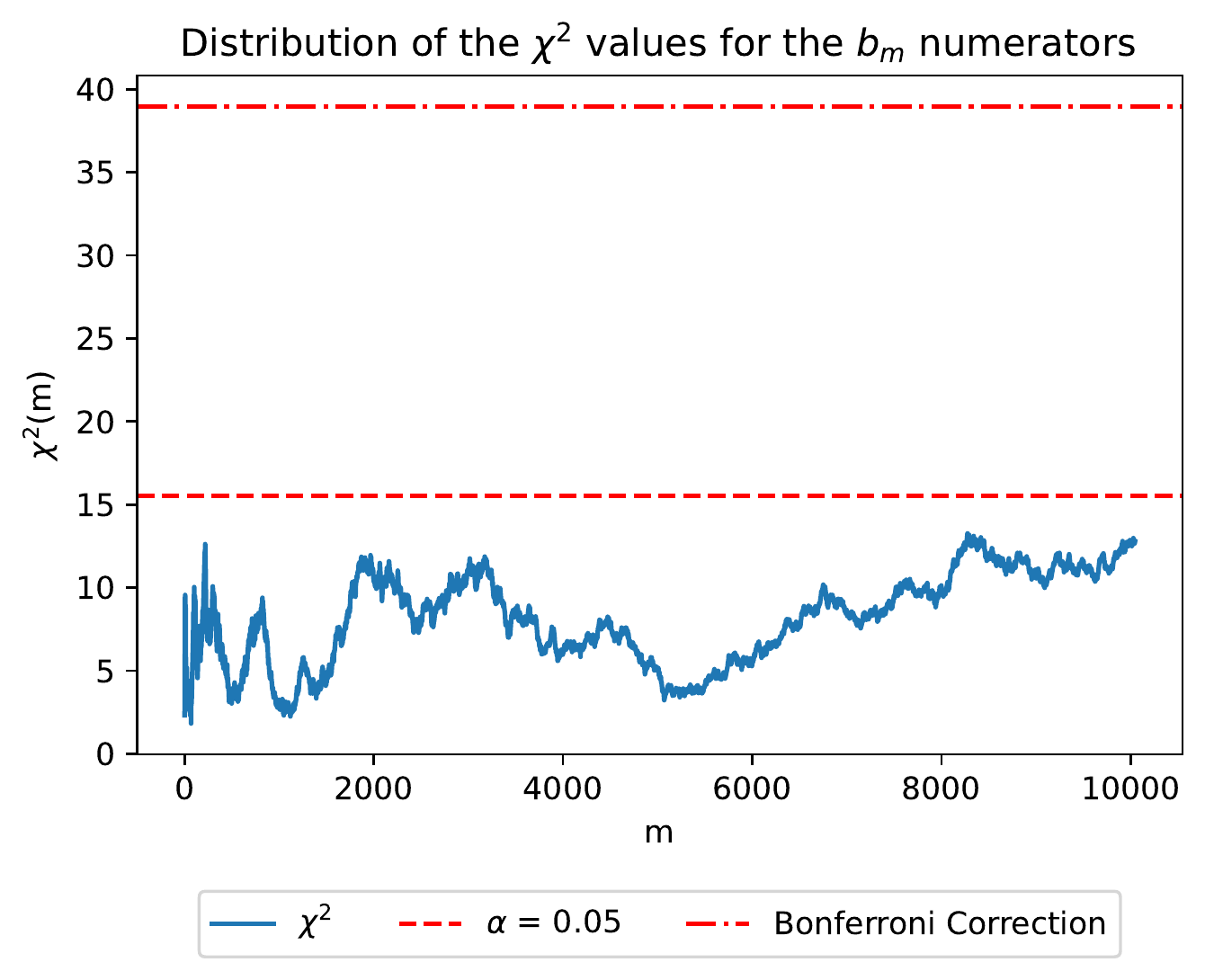} }}%
    \caption{$\chi^2$ distribution for the $a_m$ and $b_m$ numerators.}%
    \label{fig:ChiSquNum}
\end{figure}

\begin{figure}[H]%
    \subfloat[\centering $a_m$]{{\includegraphics[height = 6cm, width = 6cm]{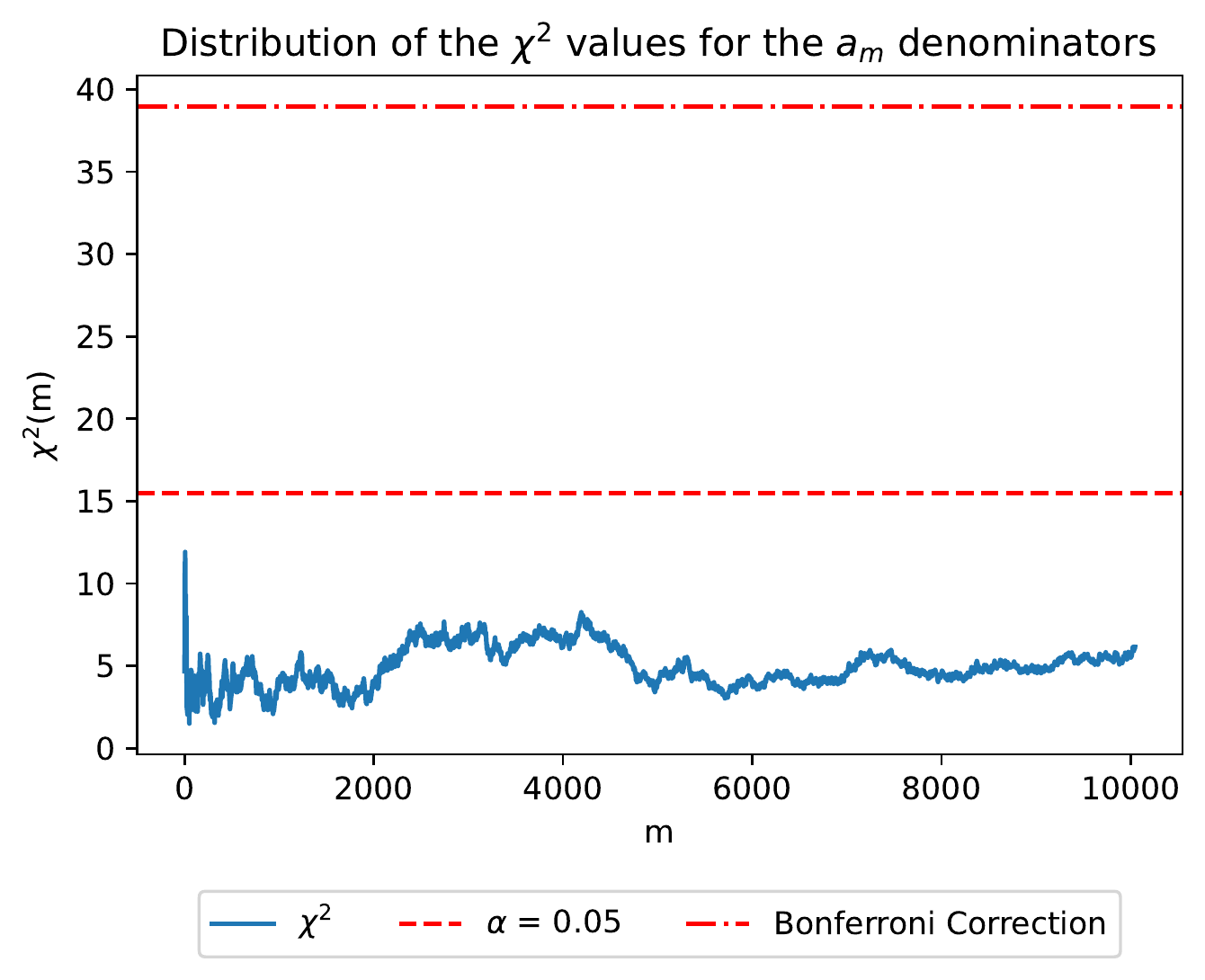} }}%
    \qquad
    \subfloat[\centering $b_m$]{{\includegraphics[height = 6cm, width = 6cm]{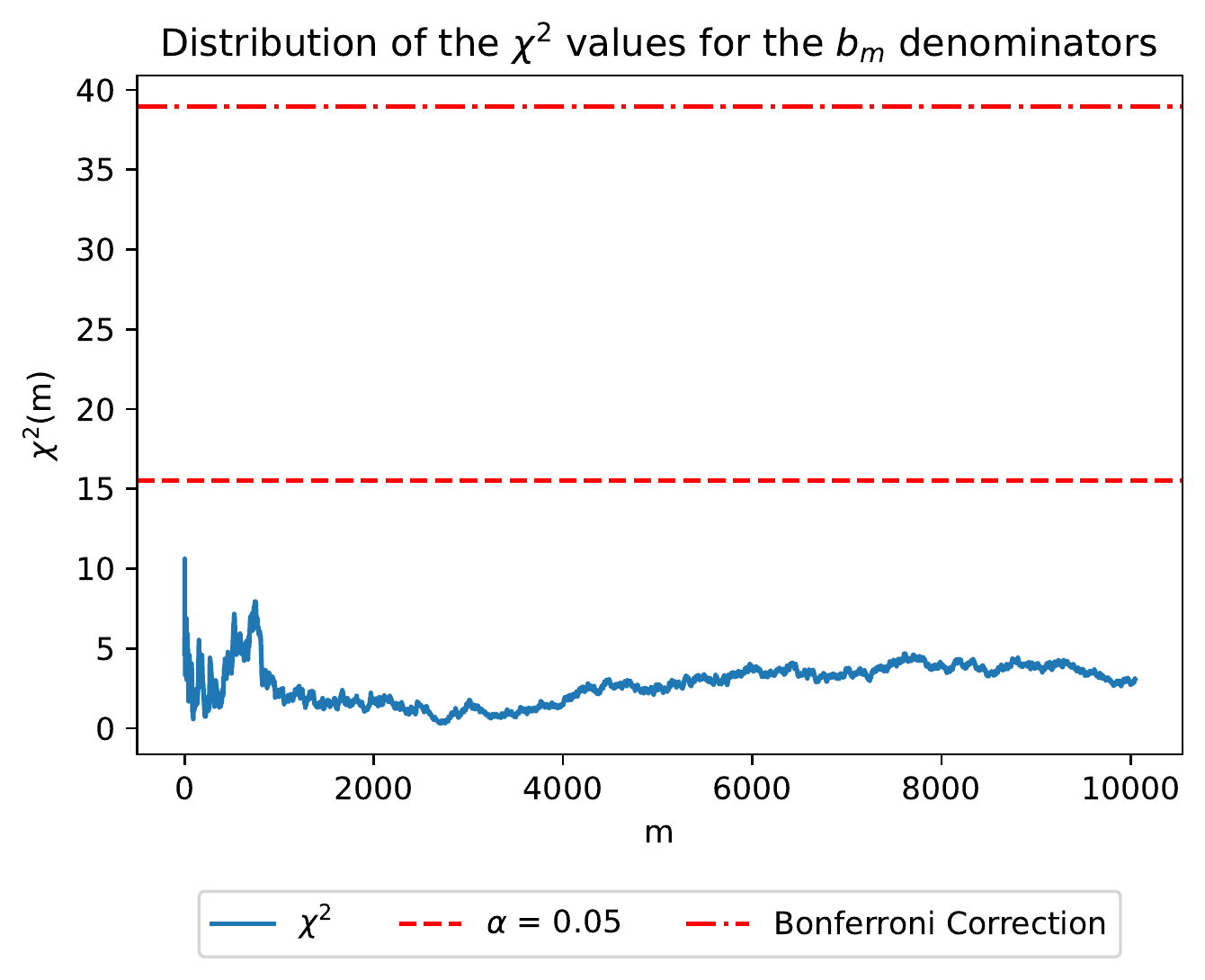} }}%
    \caption{$\chi^2$ distribution for the $a_m$ and $b_m$ denominators.}%
    \label{fig:ChiSquDen}
\end{figure}

\begin{figure}[H]%
    \subfloat[\centering $a_m$]{{\includegraphics[height = 6cm, width = 6cm]{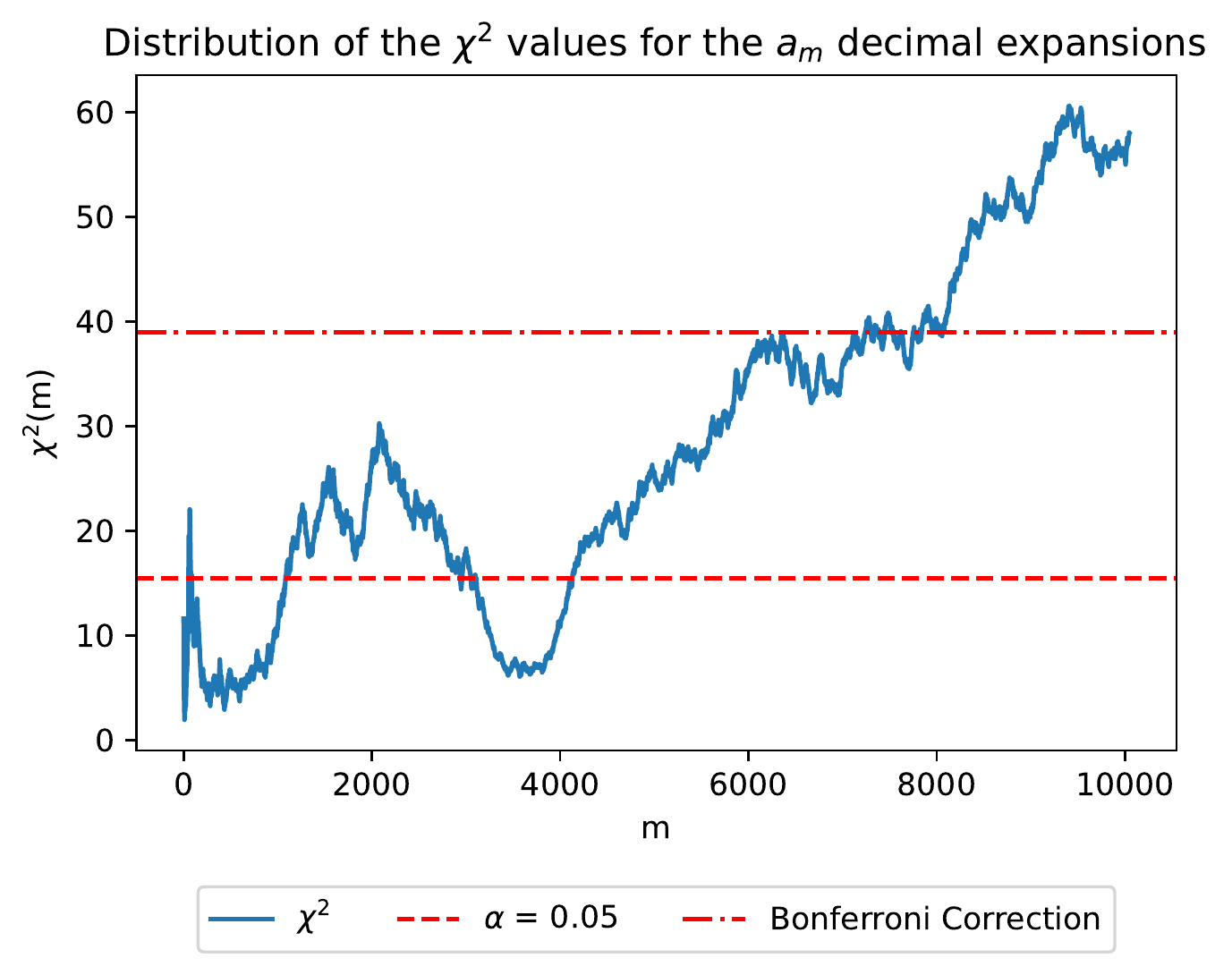} }}%
    \qquad
    \subfloat[\centering $b_m$]{{\includegraphics[height = 6cm, width = 6cm]{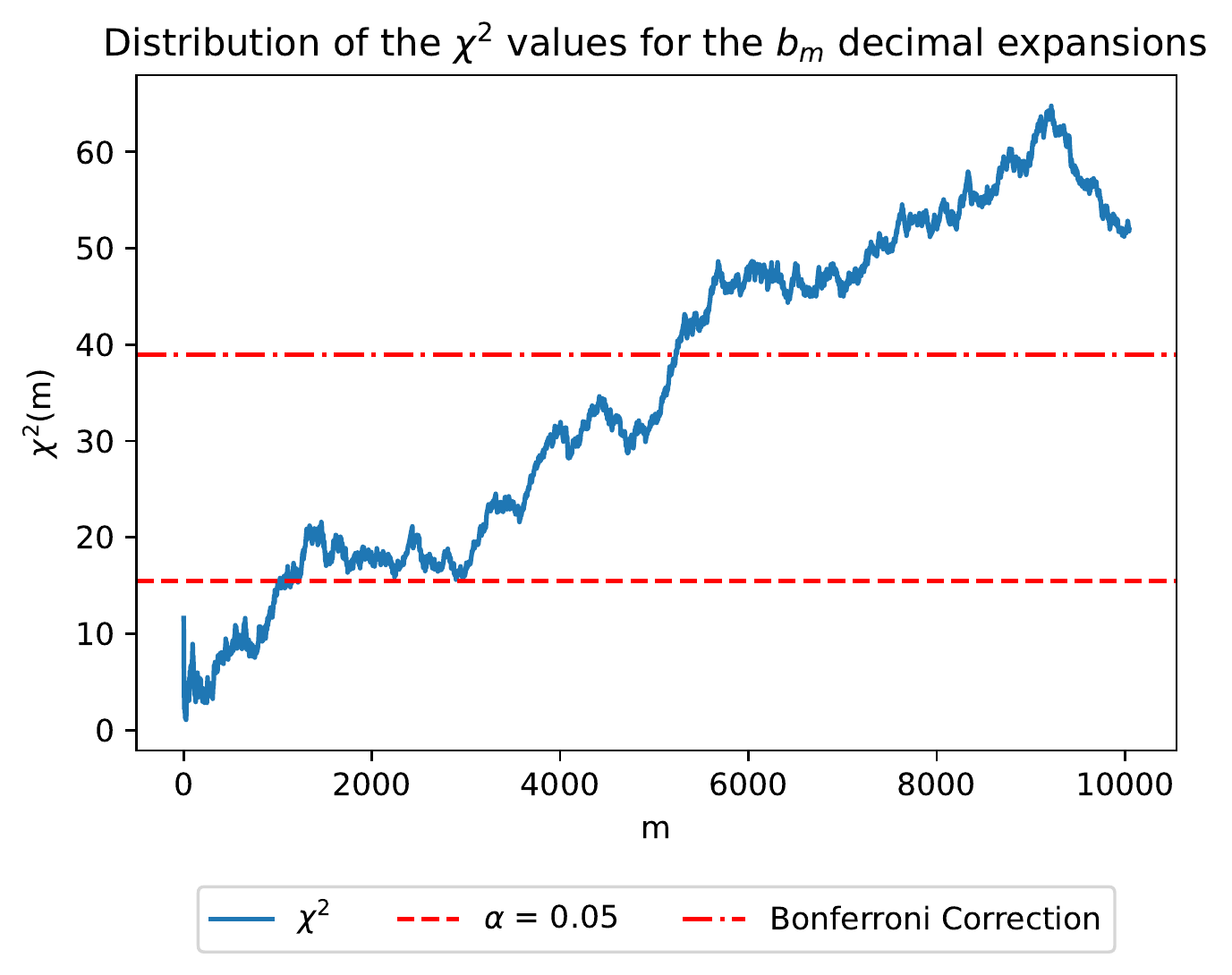} }}%
    \caption{$\chi^2$ distribution for the $a_m$ and $b_m$ decimal expansions.}%
    \label{fig:ChiSquDec}
\end{figure}

\begin{figure}[H]%
    \subfloat[\centering $a_m$]{{\includegraphics[height = 6cm, width = 6cm]{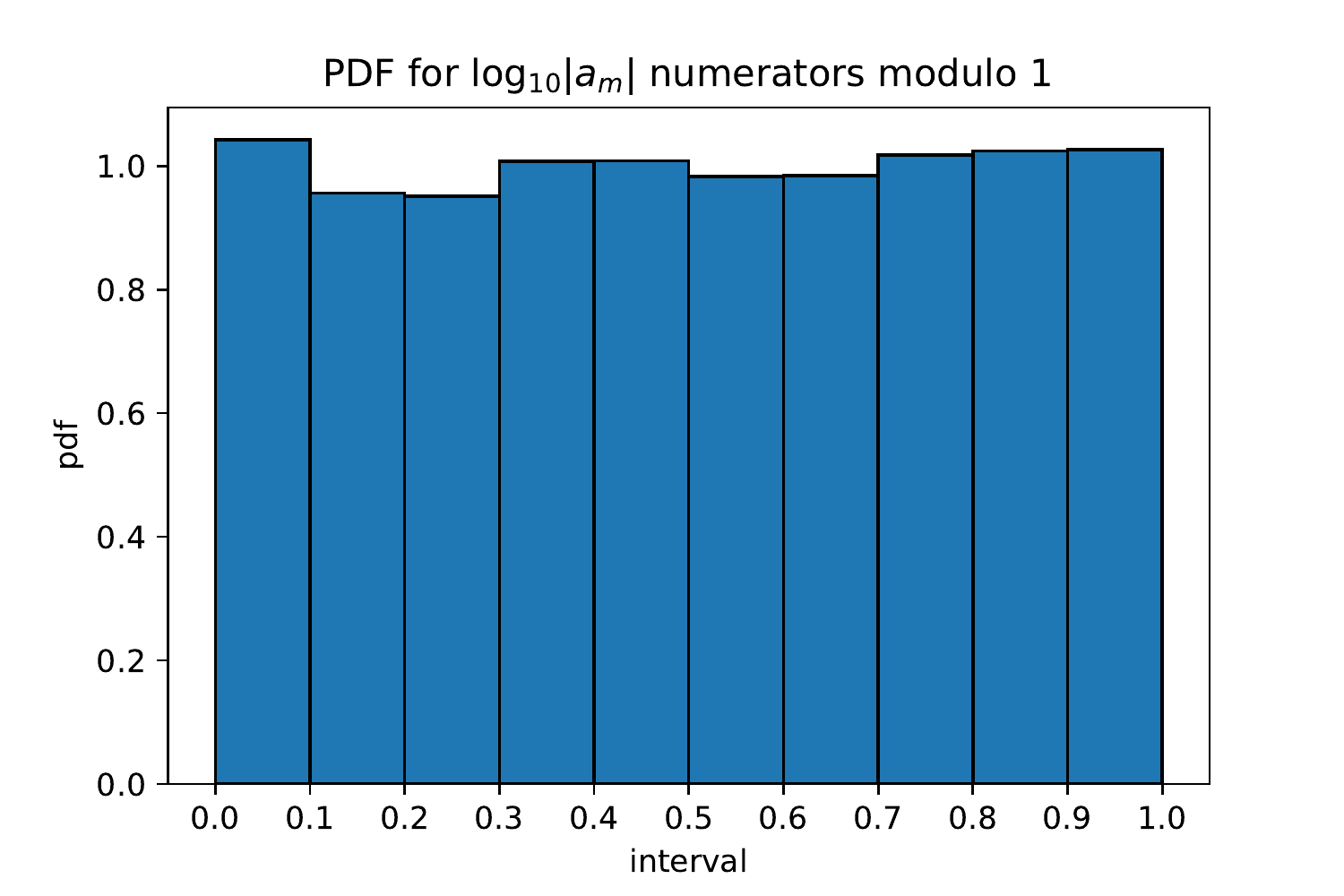} }}%
    \qquad
    \subfloat[\centering $b_m$]{{\includegraphics[height = 6cm, width = 6cm]{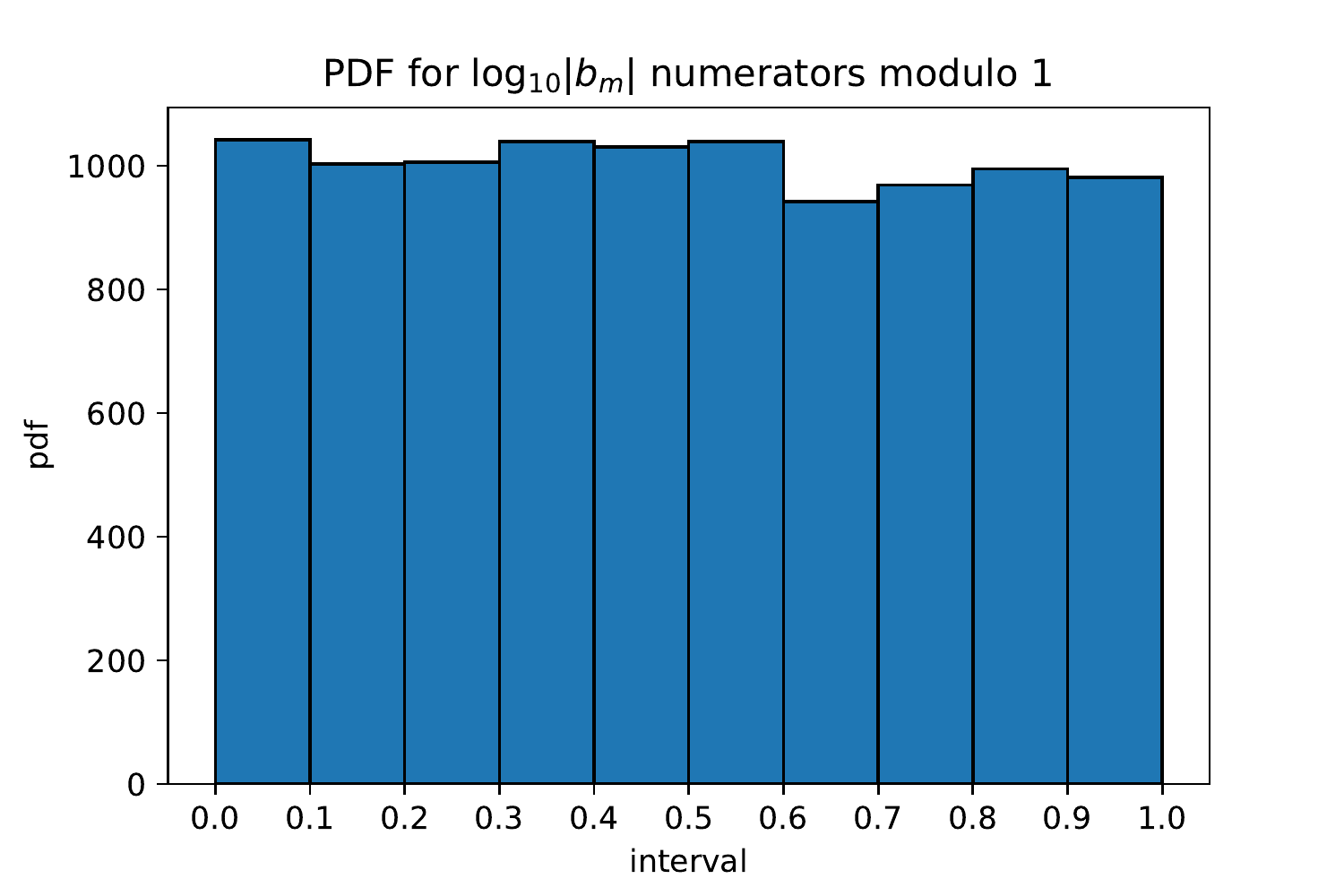} }}%
    \caption{$a_m$ and $b_m$ numerators.}%
    \label{fig:logNum}
\end{figure}

\begin{figure}[H]
    
    \subfloat[\centering $a_m$]{{\includegraphics[height = 6cm, width = 6cm]{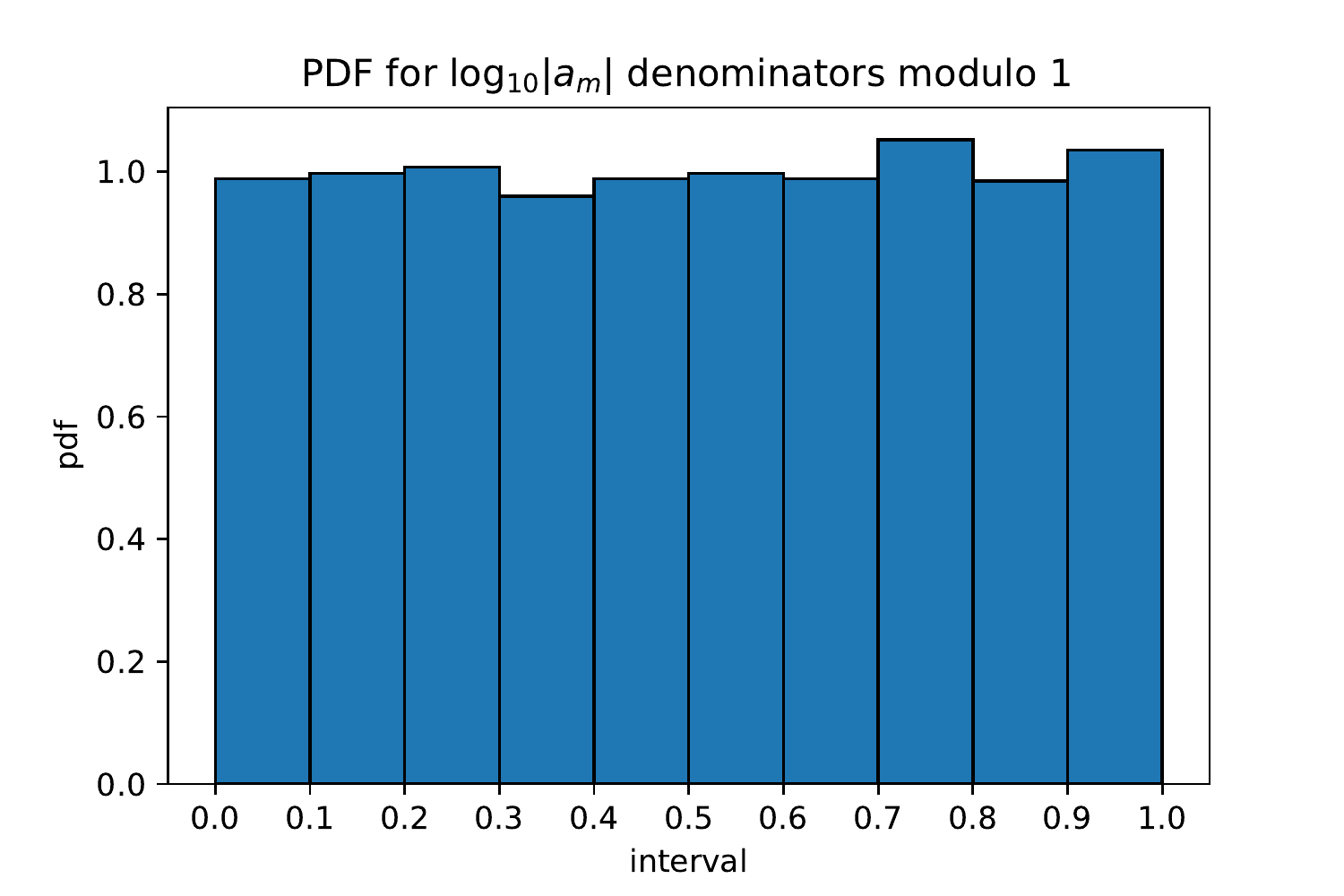} }}%
    \qquad
    \subfloat[\centering $b_m$]{{\includegraphics[height = 6cm, width = 6cm]{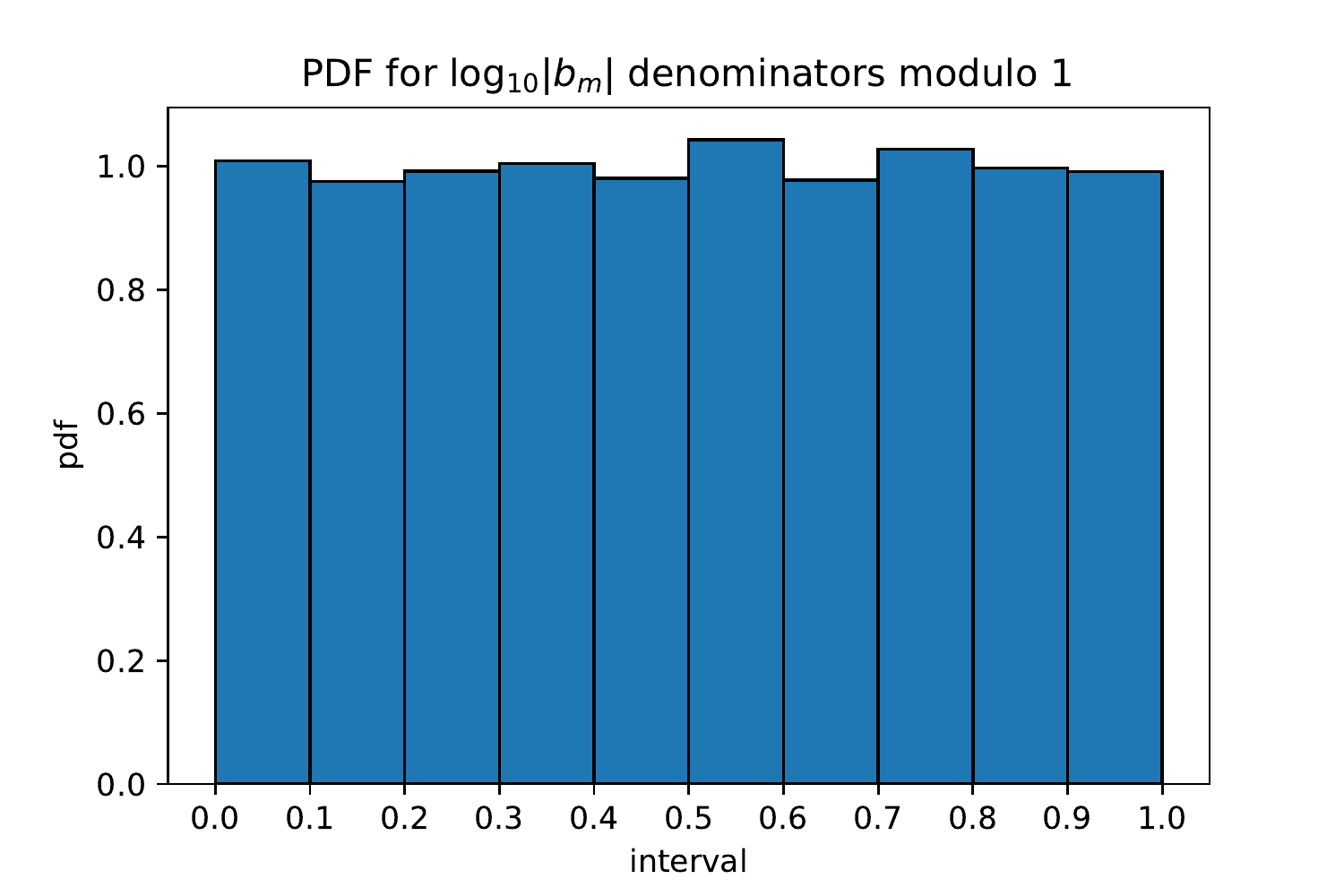} }}%
    \caption{$a_m$ and $b_m$ denominators.}%
    \label{fig:logDen}
\end{figure}

\begin{figure}[H]
    
    \subfloat[\centering $a_m$]{{\includegraphics[height = 6cm, width = 6cm]{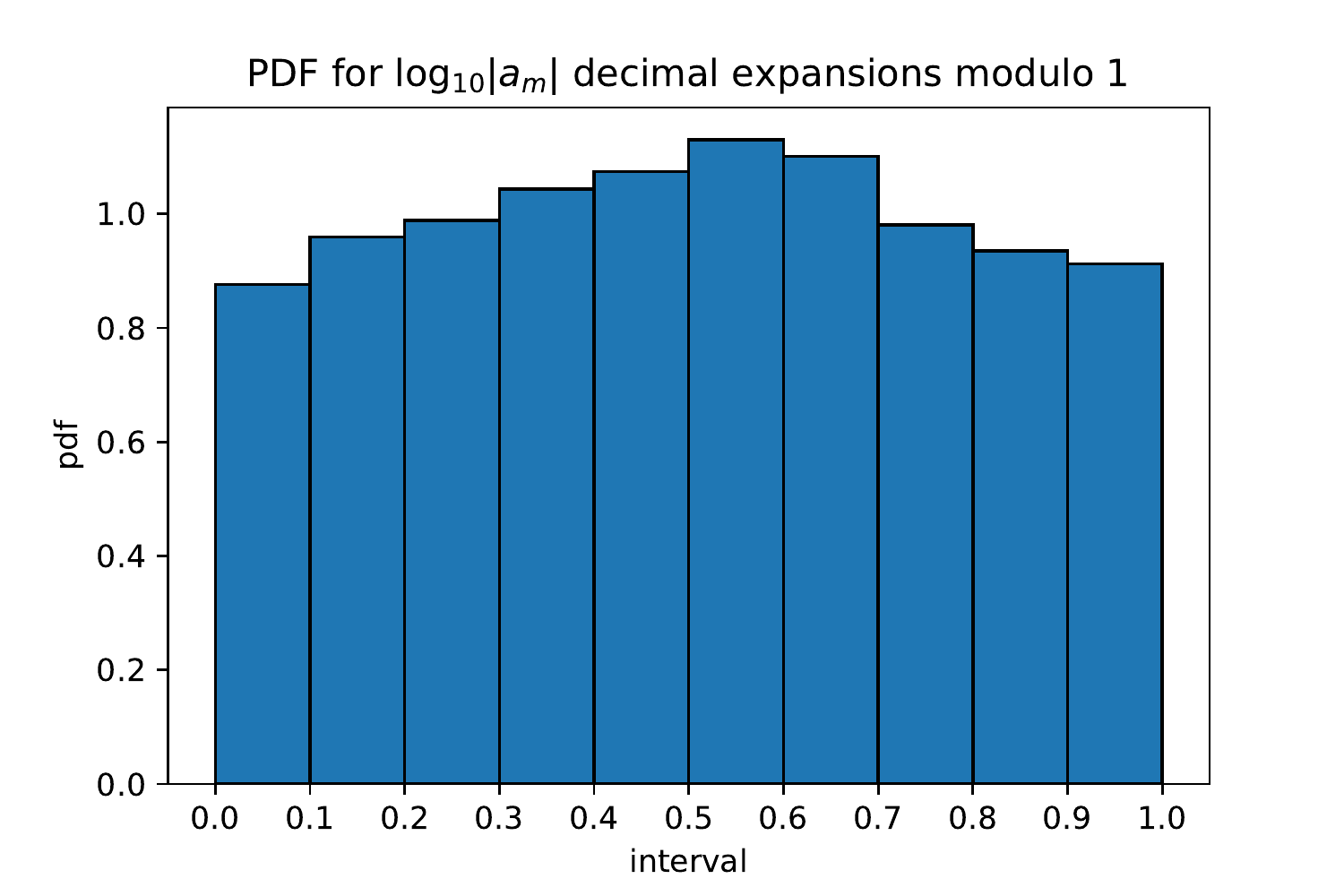}}}%
    \qquad
    \subfloat[\centering $b_m$]{{\includegraphics[height = 6cm, width = 6cm]{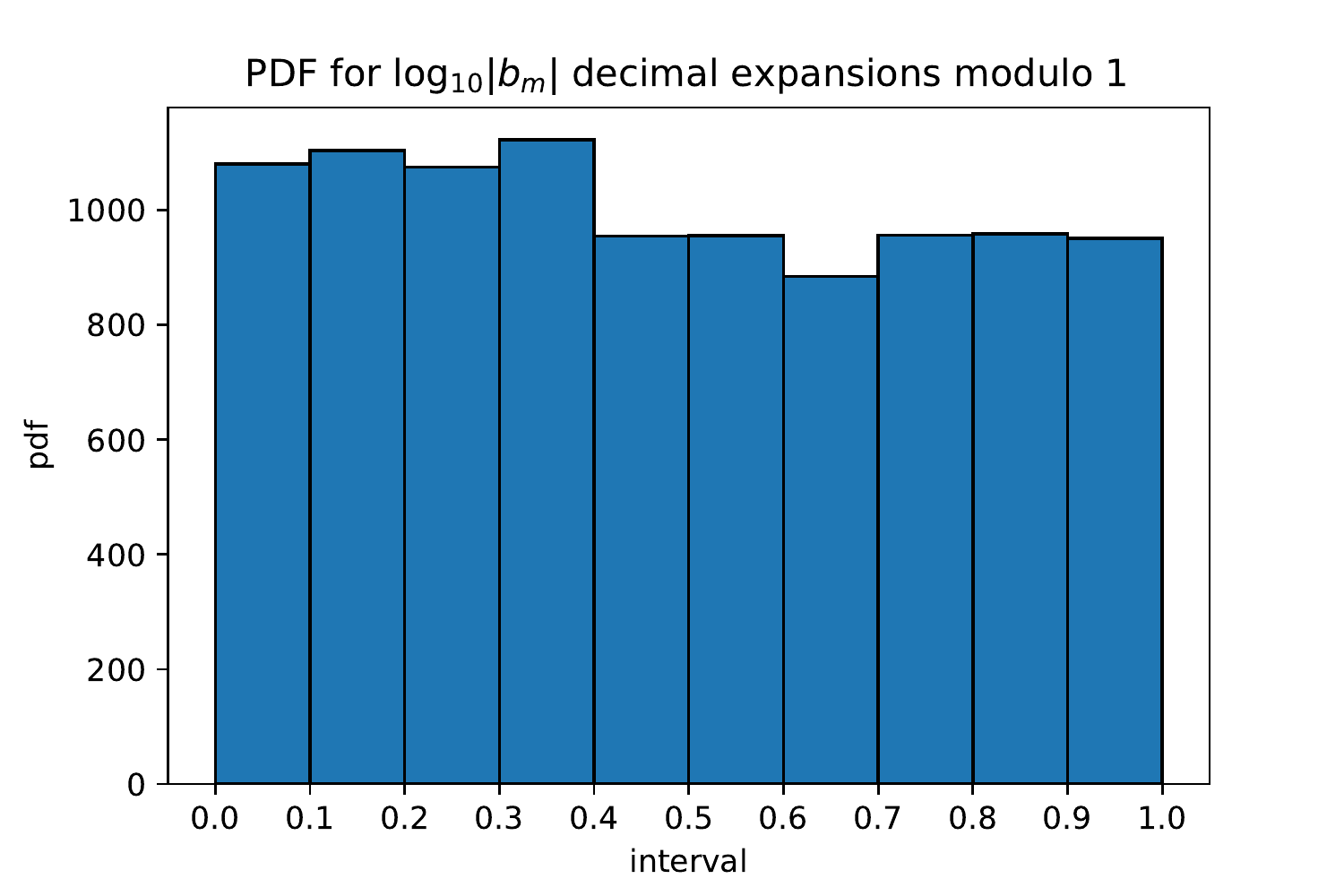}}}%
    \caption{$a_m$ and $b_m$ decimal expansions.}%
    \label{fig:logDec}
\end{figure}

We examine the distribution of the first digits of the $a_m$ and $b_m$ coefficients. As~mentioned earlier, we restrict our discussion to the non-zero coefficients. We conduct the $\chi^2$ test and distribution of the base 10 logarithms modulo 1 to~evaluate the data. The notebooks to generate the results are found under the Data Analysis folder. The~plots of the $\chi^2$ values are shown in Figures~\ref{fig:ChiSquNum}--\ref{fig:ChiSquDec}; {these figures were generated by cells 16, 17, and~18 in the Jupyter notebooks amChiSquareData.ipynb and bmChiSquareData.ipynb.}
The plots of the logarithms modulo 1 are provided in Figures~\ref{fig:logNum}--\ref{fig:logDec}; {these figures were generated by cells 10,11, and~12 in the Jupyter notebooks amLogData.ipynb and bmLogData.ipynb. Simulation results can be reproduced using the sampling\_simulation.ipynb notebook.}

{On the plots of the $\chi^2$ values, the~numerators stay below the threshold for statistical significance. The~final $\chi^2$ values for the $a_m$ and $b_m$ numerators are 5.968 and 12.785, respectively. The~type I error rates (\emph{p}-values) for the final $\chi^2$ values are 0.651 and 0.119, both of which are greater than our critical values of $\alpha = 0.05/10045$ for the individual hypotheses and $\alpha = 0.05$ for the overall test. The~sampling errors we obtain from our simulations for $a_m$ and $b_m$ are 0.043 and 0.058. The~powers relative to the null hypothesis for $a_m$ and $b_m$ are 0.355 and 0.718. As~a result, there is not sufficient evidence to reject our null hypothesis that the data are Benford.}

{
The denominators stay below the original threshold for significance. This is expected, as they consist of a random sampling of a geometric series, which is known to be Benford in most bases~\cite{MilTB1}. The~final $\chi^2$ values for the $a_m$ and $b_m$ denominators are 6.148 and 3.093, respectively. The~type I error rates (\emph{p}-values) for the final $\chi^2$ values are 0.631 and 0.928, both of which are greater than our critical values of $\alpha = 0.05/10,045$ for the individual hypotheses and $\alpha = 0.05$ for the overall test. The~sampling errors we obtain from our simulations for $a_m$ and $b_m$ are 0.047 and 0.028. The~powers relative to the null hypothesis for $a_m$ and $b_m$ are 0.366 and 0.187. As~a result, there is not sufficient evidence to reject our null hypothesis that the data are Benford.}

{
The $\chi^2$ values for the decimal expansions are the most interesting. The~$\chi^2$ values for the $a_m$ coefficients seem to increase in general, but~there are periods where the values fall. The~sequence ends above the corrected threshold for significance.  The~$b_m$ coefficients surpass the corrected threshold as well. The~final $\chi^2$ values for the $a_m$ and $b_m$ data are 58.054 and 51.934. The~type I error rates (\emph{p}-values) for the final $\chi^2$ values are $1.121$ $\times$ 10$^{-9}$ and $1.733$ $\times$10$^{-8}$, both of which are less than our critical values of $\alpha = 0.05/10,045$ for the individual hypotheses and $\alpha = 0.05$ for the overall test. The~sampling errors we obtain from our simulations for $a_m$ and $b_m$ are 0.235 and 0.243. The~powers relative to the null hypothesis for $a_m$ and $b_m$ are 0.999992 and 0.999956. As~a result, there is sufficient evidence to reject our null hypothesis that the data are Benford.}

{The distributions for the base 10 logarithms modulo 1 of the numerators are slightly skewed. There is a pattern in the $a_m$ coefficients that could account for this; we have observed that when $m = 2^n$, $a_m = 1/m$. This result seems to generalize for $a_{d,m}$, such that when $m = d^n$, $a_{d,m}$ = 1/m, which can be observed in the tables provided by Shiamuchi in~\cite{Shi2}, and~we have not found a counterexample in our computations. There is regularity in the $b_m$ numerators as discussed by Bielefeld, Fisher, and~von Haeseler in~\cite{BieFH1}, and~it is likely that similar regularities are present in the $a_m$ coefficients as well. The~$\chi^2$ values for the $a_m$ and $b_m$ data are 8.482 and 10.203, respectively. These correspond to \emph{p}-values of 0.486 and 0.334; the powers relative to the null hypothesis are 0.482 and 0.574. As a result, there is not sufficient evidence to reject the null hyopthesis that the data are uniform.}

{
The denominators consist of a sampling of integer powers of 2. Since $\log_{10}(2)$ is irrational, the~sequence $x_n = 2^n$ is Benford in base 10, and~$\log_{10}(2^n)$ (mod 1) converges to a uniform distribution~\cite{MilTB1}. Since the denominators span many orders of magnitude, it is expected that they will similarly converge in distribution. The~$\chi^2$ values for the $a_m$ and $b_m$ data are 6.334 and 4.416. These correspond to \emph{p}-values of 0.706 and 0.882; the powers relative to the null hypothesis are 0.358 and 0.248. As a result, there is not sufficient evidence to reject the null hypothesis that the data are uniform.}

It is worth noting that the distributions of the logarithms modulo 1 for $a_m$ and $b_m$ decimal expansions are skewed towards different halves of the interval. {This asymmetry may be related to how the series represent coefficients of reciprocal functions and how they may be computed from each other.}
 {The $\chi^2$ values for decimal expansions are 64.261 and 60.757. These correspond to \emph{p}-values of $2.008$ $\times$ 10$^{-10}$ and $9.580$ $\times$ 10$^{-10}$; the powers relative to the null hypothesis are 0.99998 and 0.99994. There is  sufficient evidence to reject the null hyopthesis that the data are uniform.} 

We may also investigate the magnitude of the data by computing the arithmetic mean and standard deviation of $\log_{10}{\lvert x_n \rvert}$. It is typical, but~not necessary, for~a data set to be Benford if it spans many orders of magnitude (see  Chapter 2 of~\cite{Mil1,BerH1} for an analysis that a sufficiently large spread is not enough to ensure Benfordness). Our findings for the Taylor and Laurent coefficients are summarized in Table~\ref{tab:mean,sd}. {The data are generated by cell 8 in the Jupyter notebooks amLogData.ipynb and bmLogData.ipynb, which may be found under the Data Analysis folder.}

\begin{table}[H]
\caption{Parameters for the $\log_{10}$ distribution of the data~sets.}\label{tab:mean,sd}
\newcolumntype{C}{>{\centering\arraybackslash}X}
\begin{tabularx}{\textwidth}{CCC}

\toprule
\textbf{Data Set }& \boldmath{$\mu$} & \boldmath{$\sigma$} \\ 
\midrule
$a_m$ numerators & {1899.033} & {1793.427} \\ 
\midrule
$a_m$ denominators & {1903.545} & {1793.680} \\ 
\midrule
$a_m$ decimals & {$-$4.513} & {0.679} \\ 
\midrule
$b_m$ numerators & {1899.284} & {1793.752}\\ 
\midrule
$b_m$ denominators & {1904.132} &{1793.994} \\ 
\midrule
$b_m$ decimals & {$-$4.848} & {0.689} \\ 
\bottomrule
\end{tabularx}   
%
%\vspace{15pt}
%
%\vspace{-15pt}

\end{table}

The numerators and denominators span many orders of magnitude, while the decimals do not. The~mean for the decimal expansions being negative indicates that the denominators are larger than the numerators, on average. The~ratio between the growth rates of the numerators and denominators likely has some form of regularity as well to account for the small standard deviation, but~more analysis would be needed to determine the exact relationship. These observations are consistent with the previously discussed conjecture that $0 < |b_m| < 1/m$ for all  $m$, and~it is plausible a similar relation holds for the $a_m$ coefficients as well. Ultimately, this provides insight into the growth of the coefficients and the shape of the data. Our testing provides evidence for convergence of the numerators and denominators to a Benford distribution, but there is not sufficient evidence for convergence in the decimal~expansions.

%%%%%%%%%%%%%%%%%%%%%%%%%%%%%%%%%%%%%%%%%%%%%%%%%%%%%%%%%%%%%%%%%%%%%%%%%%%%%%%%%%%%%%%%%%%%%%%%%%%%%%%%%%%%%%%%%%%%%%%%%%%%%%%%%%%%%%%%%%%%%%%%%%%
%%%%%%%%%%%%%%%%%%%%%%%%%%%%%%%%%%%%%%%%%%%%%%%%%%%%%%%%%%%%%%%%%%%%%%%%%%%%%%%%%%%%%%%%%%%%%%%%%%%%%%%%%%%%%%%%%%%%%%%%%%%%%%%%%%%%%%%%%%%%%%%%%%%
%%%%%%%%%%%%%%%%%%%%%%%%%%%%%%%%%%%%%%%%%%%%%%%%%%%%%%%%%%%%%%%%%%%%%%%%%%%%%%%%%%%%%%%%%%%%%%%%%%%%%%%%%%%%%%%%%%%%%%%%%%%%%%%%%%%%%%%%%%%%%%%%%%%
\section{On the Taylor and Laurent~Coefficients}\label{sec_coeffs}

{
This Section deals with observations and theorems on the coefficients from various authors. Our goal is to compile and highlight important results from disparate sources. We link to the papers where the original observations may be found and their proofs when applicable. Section~\ref{subsec:amRemarks} deals specifically with observations related to the $a_m$ coefficients, Section~\ref{subsec:bmRemarks} deals with the $b_m$ coefficients, and~Section~\ref{subsec:ourRemarks} deals specifically with new observations we make.}

Theorem~\ref{car} highlights the relevance of the study of the Riemann mappings $\Psi$ and $\Phi$. Much effort has been put into the understanding of the behavior of both series. We now refer to a few important results and introduce new conjectures on the behavior of the coefficients. {One of the most important theorems relating to both sets of coefficients is the following.}

\begin{Theorem} \label{padic}
The {$a_{d,m}$ and $b_{d,m}$} coefficients are $d$-adic rational numbers. In~other words, $a_{d,m}$ and $b_{d,m}$ are of the form
$$\frac{p}{d^{-v}},$$
where $p$ is an integer indivisible by $d$. The~integer $\nu$ is, by~definition, the~$d$-adic valuation $\nu_d$ of $a_{d,m}$ or $b_{d,m}$.
\end{Theorem}

{This is a combination of Theorems outlined by Shiamuchi in~\cite{Shi3,Shi4}, and~it provided our motivation for studying the numerators, denominators, and~decimal expansions, separately.} {We} consider only the case $d=2$. The~majority of the {following} results hold also for a general integer $d \ge 2$, under~simple~modifications.  

%%%%%%%%%%%%%%%%%%%%%%%%%%%%%%%%%%%%%%%%%%%%%%%%%%%%%%%%%%%%%%%%%%%%%%%%%%%%%%%%%%%%%%%
%%%%%%%%%%%%%%%%%%%%%%%%%%%%%%%%%%%%%%%%%%%%%%%%%%%%%%%%%%%%%%%%%%%%%%%%%%%%%%%%%%%%%%%
%%%%%%%%%%%%%%%%%%%%%%%%%%%%%%%%%%%%%%%%%%%%%%%%%%%%%%%%%%%%%%%%%%%%%%%%%%%%%%%%%%%%%%%
%%%%%%%%%%%%%%%%%%%%%%%%%%%%%%%%%%%%%%%%%%%%%%%%%%%%%%%%%%%%%%%%%%%%%%%%%%%%%%%%%%%%%%%
\subsection{Results for the Taylor~Coefficients}
\label{subsec:amRemarks}

\begin{Remark}
For any integers $k$ and $\nu$ satisfying $k\geq 1$ and $2^{\nu}\geq k+1$, let $m= (2k+1)2^{\nu}$. Then $a_{2,m} = 0.$
\end{Remark}

It is unknown whether the converse is true. The~proof may be found in~\cite{EwiS1}. The~authors have reported that their computation of 1000 terms of $a_{2,m}$ has not produced a zero-coefficient besides those indicated by the theorem, which is consistent with our observations. The~result may also be expressed as {the} following corollary:

\begin{Corollary} \label{cor1}
Let $m=m_{0}2^{n}$ with $n\geq 0,$ and $m_{0}$ odd. If~$3\leq m_{0}\leq 2^{n+1}$, then $a_{2,m}=0$.
\end{Corollary}

Making use of the 2-adic valuation, it is possible to obtain the following theorem, {as outlined in~\cite{Shi2}.}

\begin{Theorem} \label{thm:nuineq}
We have $-\nu_2(a_{m})\ \leq \ \nu_2((2m-2)!)$ for all $m$, with~equality attained exactly when $m$ is odd.
\end{Theorem}

{In the following, since} our interest is only for the 2-adic valuation, we will make use of the notation $\nu(x) := \nu_2(x)${, and~we immediately obtain the following remark}.
 
\begin{Remark} \label{fact}
In the case that $m$ is odd, we may obtain an efficient algorithm to compute $-\nu(a_{m})$ {through} $\nu((2m-2)!)$. {Following immediately from Theorem \ref{thm:nuineq} by }the properties of the $d$-adic evaluation outlined in~\cite{Shi3}, we have,
\begin{align}
    \nu((2m-2)!) \ = \ \sum_{l=1}^\infty \left \lfloor{\frac{2m-2}{2^l}}\right \rfloor \ = \ \sum_{l=0}^\infty \left \lfloor{\frac{m-1}{2^l}}\right \rfloor.
\end{align}

Therefore, if~we set a value, $N$, the~denominator's exponent for every odd number $m < 2^N$ is given by
\begin{align}
    -\nu(a_{m}) \ = \ \nu((2m-2)!) \ = \ \sum_{l=0}^N \left \lfloor{\frac{m-1}{2^l}}\right \rfloor.
\end{align}
\end{Remark}

{We may also summarize} Theorem 3.1 and Corollary  3.5 from~\cite{Shi4} for the case that $d=2$.

\begin{Theorem}
 Given $m \in \mathbb{N}\sm\{1\}$, we have that $-\nu(a_m) \le x(m)$, where
 \begin{align*}
     x(m)  \ = \  \nu((m-1)!) + m-1.
 \end{align*}
 
 Under the same assumptions, the~result is also true with
 \begin{align*}
     x(m) \ = \ \left \lceil \nu(m-1) + m - 1 \right \rceil.
 \end{align*}
\end{Theorem}

%%%%%%%%%%%%%%%%%%%%%%%%%%%%%%%%%%%%%%%%%%%%%%%%%%%%%%%%%%%%%%%%%%%%%%%%%%%%%%%%%%%%%%%
%%%%%%%%%%%%%%%%%%%%%%%%%%%%%%%%%%%%%%%%%%%%%%%%%%%%%%%%%%%%%%%%%%%%%%%%%%%%%%%%%%%%%%%
%%%%%%%%%%%%%%%%%%%%%%%%%%%%%%%%%%%%%%%%%%%%%%%%%%%%%%%%%%%%%%%%%%%%%%%%%%%%%%%%%%%%%%%
%%%%%%%%%%%%%%%%%%%%%%%%%%%%%%%%%%%%%%%%%%%%%%%%%%%%%%%%%%%%%%%%%%%%%%%%%%%%%%%%%%%%%%%
\subsection{Results for the Laurent~Coefficients} \label{bm} 
\label{subsec:bmRemarks}

Similar results hold for the $b_m$ coefficients. We use the notation $m = 2^n m_0$, where $m_0$ is odd. The~first result was presented in~\cite{Lev} in~1988.

\begin{Remark} \label{cor2}
 If $n\geq 2$ and $m_{0} \leq 2^{n+1}-5$, then $b_{m} = 0$.
\end{Remark}

It is still unknown whether the converse of this theorem is true. In~\cite{EwiS2}, the~only coefficients that have been observed to be zero are those mentioned in this theorem. The~following result, from~\cite{BieFH1}, underlines that a result similar to the one for the $a_m$ holds.

\begin{Theorem} \label{sub1}
We have $-\nu(b_{m}) \leq \nu((2m+2)!)$ for all $m$, and~equality is attained exactly when $m$ is odd.
\end{Theorem}

Don Zagier has made several observations and conjectures about the exponents of $b_m$. {We shall later} extend them to the $a_m$ coefficients. The~original conjectures are outlined in~\cite{BieFH1}.

\begin{Conjecture} \label{con:bmperiod}
For the $b_m$ coefficients, we have the~following.
\begin{itemize}
    \item If $n \ = \ 0$, then $-\nu(b_m) \ = \ \nu((2m+2)!)$.
    \item If $n \  = \ 1$, then $-\nu(b_m) \ = \ \nu(((2m+2)/3)!) + \epsilon(m_0)$, where $\epsilon(m_0) \ = \ 0$ if $m_0 \ \equiv \ 11$ mod 12 and $1$ otherwise.
    \item If $n  \ =  \ 2$, then $-\nu(b_m) \ = \ \nu(((2m-25)/7)!) + \epsilon(m_0)$, where $\epsilon(m_0)$ moves with periodicity of 28, as~expressed in~\cite{BieFH1}.
\end{itemize}
\end{Conjecture}

%%%%%%%%%%%%%%%%%%%%%%%%%%%%%%%%%%%%%%%%%%%%%%%%%%%%%%%%%%%%%%%%%%%%%%%%%%%%%%%%%%%%%%%
%%%%%%%%%%%%%%%%%%%%%%%%%%%%%%%%%%%%%%%%%%%%%%%%%%%%%%%%%%%%%%%%%%%%%%%%%%%%%%%%%%%%%%%
%%%%%%%%%%%%%%%%%%%%%%%%%%%%%%%%%%%%%%%%%%%%%%%%%%%%%%%%%%%%%%%%%%%%%%%%%%%%%%%%%%%%%%%
%%%%%%%%%%%%%%%%%%%%%%%%%%%%%%%%%%%%%%%%%%%%%%%%%%%%%%%%%%%%%%%%%%%%%%%%%%%%%%%%%%%%%%%
\subsection{New Remarks on the~Coefficients}

\label{subsec:ourRemarks}

\begin{Remark} \label{subseq}
{Theorem}~\ref{sub1} and the results on $\nu(x!)$ reported in Remark~\ref{fact} imply that the denominator's coefficients for $b_m$ have regularity in terms of consecutive differences for odd $m$. Zagier's conjectures posit similar patterns for even $m$. This gives a constant discrete derivative in the $b_m$ denominator's exponents for $m \equiv 2 \text{ mod 4} \; (n=1)$, $m \equiv 4 \text{ mod 8} \; (n=2)$, $m \equiv 8 \text{ mod 16}\; (n=3); \ \dots \ m \equiv 2^N \text{mod} \ 2^{N+1} (n = N)$. Note that for each of these subsequences, there are some starting values of $0$. These correspond exactly with the values predicted by Corollaries~\ref{cor1} and~\ref{cor2}, respectively; when $a_m = 0$, the~algorithm gives $0$ for the denominator's exponent.
In general, for~each $n$ there is a partial periodicity with period $2(2^{n+1} - 1)$ in $m_0$, and~equivalently, $2^{n+1}(2^{n+1} - 1)$ in $m$. 
\end{Remark}

Another direction is to calculate the slope of each of the subsequences, since they seem to grow linearly. Our observations make use of the previous remark and have led to  the following~conjecture.

\begin{Conjecture} \textbf{\ref{slope}}
{Let m be written as $2^nm_0$ as before, with~$n=\overline{n}$ fixed. Then, the sequence $\{-\nu(a_m)\}_{n=\overline{n}}$ is asymptotically linear, with~slope $2/(2^{\overline{n}+1}-1).$}
\end{Conjecture}

The numerators tend to follow similar behavior as the denominators. In~particular, the~modulus of the numerators in the subsequences tend to organize as follows: $\{a_m\}_{n=0}$ > $\{a_m\}_{n=1}$ > $\cdots$  $\{a_m\}_{n=N}$. The~possibility of bounding the numerators by making use of its associated denominator is a subject of further~study.

We now {extend Conjecture \ref{con:bmperiod} to the $a_m$ coefficients, as follows:}

\begin{Conjecture}
For the $a_m$ coefficients we~have:
\begin{itemize}
    \item if $n \ = \ 0$, $-\nu(a_m)  \ = \ \nu((2m-2)!)$, and~    \item if $n  \ = \ 1$, $-\nu(a_m)  \ = \  \nu(((2m-2)/3)!) + \epsilon(m_0)$, where $\epsilon(m_0) = 1$ if $m_0 \not \equiv3$ mod 12. Otherwise, it follows the pattern in Table~\ref{Tab1}.
\end{itemize}
\end{Conjecture}

This suggests a partial periodicity with period $2^4(2^{n+1} - 1)$ in $m_0$, or~of $2^{n+4}(2^{n+1} - 1)$ in $m$. 
As before, it is possible to write $m_0$ as $2(2^{n+1} - 1)k + l$, but~it is more difficult to identify a general pattern in this~case.

\begin{table}[H]
\caption{The distribution of $\epsilon(m_0)$ when $m_0 \equiv_{12}3$ has periodicity $16$ $\times$ 12$=192$. From~$m_0 = 195$, it repeats itself, and~will have the same $\epsilon$($m_0$) of $m_0=3$ and~following.}
\label{Tab1}
\newcolumntype{C}{>{\centering\arraybackslash}X}
\begin{tabularx}{\textwidth}{CC}
\toprule
 \boldmath{$m_0$ \text{(\textbf{mod 192})}} & \boldmath{$\epsilon$($m_0$)} \\
\midrule
 3  &  $-$2  \\
 15 &  $-$2    \\
 27 &  $-$4   \\
 39 &  $-$3   \\
 51 &  $-$2  \\
 63 &  $-$2   \\
 75 &  $-$2    \\
 87 &  $-$3    \\
 99 &  $-$4      \\
 111 &  $-$2    \\
 123 &  $-$2    \\
 135 &  $-$6        \\
 147 &  $-$3\\
 159 &  $-$2\\
 171 &   $-$2\\
 183 &  $-$3  \\
\bottomrule
\end{tabularx}
\end{table}

%%%%%%%%%%%%%%%%%%%%%%%%%%%%%%%%%%%%%%%%%%%%%%%%%%%%%%%%%%%%%%%%%%%%%%%%%%%%%%%%%%%%%%%%%%%%%%%%%%%%%%%%%%%%%%%%%%%%%%%%%%%%%%%%%%%%%%%%%%%%%%%%%%%
%%%%%%%%%%%%%%%%%%%%%%%%%%%%%%%%%%%%%%%%%%%%%%%%%%%%%%%%%%%%%%%%%%%%%%%%%%%%%%%%%%%%%%%%%%%%%%%%%%%%%%%%%%%%%%%%%%%%%%%%%%%%%%%%%%%%%%%%%%%%%%%%%%%
%%%%%%%%%%%%%%%%%%%%%%%%%%%%%%%%%%%%%%%%%%%%%%%%%%%%%%%%%%%%%%%%%%%%%%%%%%%%%%%%%%%%%%%%%%%%%%%%%%%%%%%%%%%%%%%%%%%%%%%%%%%%%%%%%%%%%%%%%%%%%%%%%%%
\section{Future~Work}

The most natural extension of our work would be to generate more coefficients, which would allow more thorough statistical testing. Given more data, we could look at the coefficients over a certain zoom or average the coefficients over certain subsets. It would be particularly interesting if certain subsets of the coefficients also converge to {the} Benford distribution. {We may also look at the powers of the denominators to observe whether they follow a Benford distribution.}

The algorithms for computing the coefficients of the Mandelbrot can also be easily generalized to obtain other abstract Multibrot sets, which could be analyzed using the same methods. We could look at the data {in} different bases to observe whether Benfordness holds there. It would be interesting to see the numerators and decimal expansions of the coefficients for the Multibrot set of degree $d$ follow a Benford distribution in base $d$; the denominators will not since they are sampled from a geometric series with a common ratio $d$, and they are not Benford in the base of the common ratio. 

The results of Section~\ref{sec_coeffs} also present interesting extensions for future work. In~particular, Remark~\ref{subseq} suggests that dividing the coefficients {into} subsequences to be bounded separately may be the best approach to study the convergence of the Laurent series of the coefficients. This approach, which has not been followed in the past, to~the best of our knowledge, could provide valuable results in the study of the local connectedness of $\mathcal{M}$.

\vspace{6pt} 

%%%%%%%%%%%%%%%%%%%%%%%%%%%%%%%%%%%%%%%%%%
%% optional
%\supplementary{The following supporting information can be downloaded at:  \linksupplementary{s1}, Figure S1: title; Table S1: title; Video S1: title.}


%%%%%%%%%%%%%%%%%%%%%%%%%%%%%%%%%%%%%%%%%%
\authorcontributions{This paper is a result of a collaborative effort of all the authors on all aspects of the work, with~J.D. on visualization and formal analysis, W.F. on data curation and software, and~F.B. on formal analysis and conceptualization. The~research was performed as part of the 2021 Polymath Jr program, advised by T.C.M., S.J.M. and D.S. All authors have read and agreed to the published version of the manuscript.}

\funding{This work was partially supported by NSF Grant~DMS2113535.}

\institutionalreview{Not applicable.}

\informedconsent{Not applicable.}

\dataavailability{The codes and data can be found here \url{https://github.com/DannyStoll1/polymath-fractal-geometry} (accessed on 7 February 2023).}

\acknowledgments{We thank our colleagues from the 2021 Polymath Jr REU program for helpful comments on earlier~drafts.}

\conflictsofinterest{The authors declare no conflicts of~interest.}

%%%%%%%%%%%%%%%%%%%%%%%%%%%%%%%%%%%%%%%%%%
\begin{adjustwidth}{-\extralength}{0cm}

\reftitle{References}

\end{adjustwidth}
\end{document}